\documentclass[english,preprint,3p,sort&compress]{elsarticle}
\usepackage[T1]{fontenc}
\usepackage[latin9]{inputenc}
\usepackage{bm}
\usepackage{amsmath}
\usepackage{graphicx}
\usepackage{esint}

\makeatletter

\newcommand{\lyxdot}{.}

\usepackage{lineno,hyperref}
\modulolinenumbers[5]

\makeatother

\usepackage{babel}
\begin{document}

\begin{frontmatter}{}

\title{Parallel optimized sampling for stochastic equations}

\author{Bogdan Opanchuk, Simon Kiesewetter and Peter D. Drummond }

\address{Centre for Quantum and Optical Science, Swinburne University of Technology,
Melbourne, Australia.}

\begin{abstract}
Stochastic equations play an important role in computational science,
due to their ability to treat a wide variety of complex statistical
problems. However, current algorithms are strongly limited by their
sampling variance, which scales proportionate to $1/N_{S}$ for $N_{S}$
samples. In this paper, we obtain a new class of variance reduction
methods for treating stochastic equations, called parallel optimized
sampling. The objective of parallel optimized sampling is to reduce
the sampling variance in the observables of an ensemble of stochastic
trajectories. This is achieved through calculating a finite set of
observables~\textemdash{} typically statistical moments~\textemdash{}
in parallel, and minimizing the errors compared to known values. The
algorithm is both numerically efficient and unbiased. Importantly,
it does not increase the errors in higher order moments, and generally
reduces such errors as well. The same procedure is applied both to
initial ensembles and to changes in a finite time-step. Results of
these methods show that errors in initially optimized moments can
be reduced to the machine precision level, typically around $10^{-16}$
in current hardware. For nonlinear stochastic equations, sampled moment
errors during time-evolution are larger than this, due to error propagation
effects. Even so, we provide evidence for error reductions of up to
two orders of magnitude in a nonlinear equation example, for low order
moments, which is a large practical benefit. The sampling variance
typically scales as $1/N_{S}$, but with the advantage of a very much
smaller prefactor than for standard, non-optimized methods. \end{abstract}


\end{frontmatter}{}

\section{Introduction}

Stochastic differential equations (SDEs) play a universal and important
role in many disciplines requiring quantitative modeling~\cite{Langevin1908-theorie,Lemons1997-langevin,Uhlenbeck1930-theory,Arnold1992-stochastic,Risken1996,Gardiner1997,Karatzas1991-brownian,VanKampen2007-stochastic,Glasserman2010-monte}.
Their practical advantage is that, when used to solve large statistical
problems, the ability to randomly sample greatly reduces the complexity
of treating the full distribution function. They are employed for
treating problems ranging from statistical physics, chemistry and
engineering through to economics, biology and financial modeling.
As a result of this large field of applications, there is substantial
literature on the algorithms used to solve them~\cite{Mil'shtein1975,Ermak1980-numerical,Talay2007-expansion,Kloeden1992,Burrage2000-note,Higham2001-algorithmic},
and these are generally different to methods for ordinary differential
equations (ODEs).

This utility is not without a price. Such numerical algorithms typically
utilize many independent trajectories. Any statistical result will
therefore have a sampling variance that scales $1/N_{S}$ for $N_{S}$
random samples. This causes to a typical error $\epsilon_{S}$ that
scales as $1/\sqrt{N_{S}}$, giving an error that only decreases slowly
with the total computation time. Numerical algorithms generally aim
to reduce the truncation error $\epsilon_{T}$ due to the discretization
in time, giving an error of order $\left(\Delta t\right)^{-p}$ for
methods of order $p$. Yet reducing the truncation error cannot reduce
the total error to less than the sampling error. Since the sampling
error is often a large part of the total error one has to deal with,
reducing the truncation error has a limited effectiveness.

Here we propose a new class of reduced variance SDE algorithms called
\emph{parallel} \emph{optimized sampling} or POS algorithms for short.
This approach can provide extremely useful variance reduction compared
to independent random sampling methods. The method locally eliminates
sampling errors over a finite set of moments, up to numerical precision.
Higher order moment errors are reduced also. Global sampling error
improvement over a finite time interval is not as large as this, due
to error propagation effects, but is still very substantial. Other
variance reduction techniques that can also improve sampling errors
include the use of low discrepancy (or quasi-random) sequences~\cite{Morokoff1995-quasi},
and an extrapolated hierarchy of simulations~\cite{Giles2008-multilevel},
which we do not treat here.

To achieve our results, we introduce a modification of the standard
independent stochastic trajectory approach. A set of $N_{S}$ trajectories
is solved \emph{in parallel} using an optimized algorithm. This gives
an improvement that depends on the order of the statistical moment
being calculated. The present method is a highly practical one, since
the run time is still only linear in $N_{S}$. This means that the
error reductions are of real value. If the run time scaled as $N_{S}^{2}$,
and errors were reduced to $1/N_{S}$, the reductions would not be
so useful, since the total computational overhead at fixed error would
not change.

There is a close relationship between the present approach and methods
of moment hierarchies~\cite{Risken1996} or cumulant expansions.
These existing cumulant based methods have a major drawback. One must
close the infinite hierarchy with an arbitrary ansatz at some finite
order. If this is incorrect, the method will fail~\cite{Morillo2014-checking},
and several examples exist for this problem. In fact, such cumulant
hierarchies cannot in general be truncated consistently~\cite{Hanggi1980-remark}.
By contrast, our method requires no ansatz; it gracefully reduces
in a consistent way to the usual, unbiased stochastic method at large
sample sizes. It is a unification of moment hierarchy and stochastic
methods.

We focus mostly on simple one-dimensional cases here for simplicity,
except for one higher-dimensional example to illustrate the approach.
Extensions of the POS method to cases of higher-order convergence
in time will be treated elsewhere. However, as one example of potential
applications, sampling errors can play a major role in quantum simulations
using phase-space representations. The difficulty is that distribution
tails may only be weakly bounded~\cite{Gilchrist1997-positive}.
In such cases, sampling error can determine whether a quantum simulation
is feasible or not. One approach to solving this problem is using
weighted trajectory methods~\cite{Deuar2002}. Another method, which
was successfully used in treating long-range order in the fermionic
Hubble model, is to impose global conservation laws~\cite{Aimi2007-gaussian}.
The present approach has a more general applicability than either
of these, but can also be combined with such earlier techniques.

The organization of the paper is as follows. In Section~\ref{sec:stochastic-integration-error},
we analyze the stochastic errors that can arise in computational methods
for stochastic differential equations. In Section~\ref{sec:static:initial-conditions}
sampling errors are reduced in initial observables by introducing
a new type of stochastic noise generator. Numerical examples of this
are treated in Section~\ref{sec:parallel-optimized-sampling}, showing
that variance in selected moments can be reduced to machine precision.
Next, in Section~\ref{sec:dynamic}, this approach is extended to
dynamical optimization of stochastic equation algorithms, with numerical
examples given in the remaining sections. In Section~\ref{sec:Synthetic-one-step-benchmarks}
we start by treating single-step optimization, which ignores error-propagation.
In Section~\ref{sec:benchmarks}, we demonstrate the global, or multi-step
performance of these methods for a number of one-dimensional linear
and nonlinear cases, showing that global sampling error variance is
substantially reduced compared to traditional independent trajectory
approaches, typically by at least an order of magnitude. Finally,
Section~\ref{sec:Two-dimensional-case} treats a two-dimensional
example which has very similar behaviour. Section~\ref{sec:conclusions}
gives our summary and conclusions.

\section{Stochastic integration error\label{sec:stochastic-integration-error}}

Stochastic integration involves both truncation errors due to finite
step-size, and errors due to finite samples. In this section we analyze
the performance of independent SDE methods by considering the overall
resources needed to solve the equations for some given error. To explain
the background to our approach, we show that increasing the convergence
order is a less effective way to reduce errors than one might expect,
because of a trade-off between the resources needed to reduce both
truncation and sampling errors. As large ensembles are required to
reduce sampling errors, we find that higher-order methods for SDEs,
while undoubtedly useful, are not as effective as in ODE integration.

In the rest of this paper, we show how to overcome this problem by
using parallel algorithms that reduce sampling errors.

\subsection{Distributions and observables }

We start by establishing our notation. The observables of a probability
distribution function $P\left(t,\mathbf{x}\right)$ are the fundamental
objects of interest in stochastic calculations. These are probability-weighted
integrals of functions $o_{m}\left(\mathbf{x}\right)$ over \textbf{$\mathbf{x}$}:
\begin{equation}
\left\langle o_{m}\right\rangle =\int\mathrm{d}\mathbf{x}o_{m}\left(\mathbf{x}\right)P\left(t,\mathbf{x}\right).
\end{equation}
In general, $\mathrm{d}\mathbf{x}$ is a Euclidean measure over a
$d$-dimensional real or complex space of variables $\mathbf{x}$.
The numerical examples given here will focus mainly on moments of
one-dimensional distributions with $o_{m}(x)=x^{m}$, giving a vector
whose error should be minimized:
\begin{equation}
\mathbf{o}=\left\{ o_{1},\ldots,o_{M}\right\} .
\end{equation}

With sampling methods, the time-evolution is solved by generating
a number $N_{S}$ of independent sample paths $\mathbf{x}^{(n)}(t)$.
This procedure is equivalent to approximating the true distribution
by the sampled distribution defined as
\begin{equation}
P_{S}\left(t,\mathbf{x}\right)=\frac{1}{N_{S}}\sum_{n=1}^{N_{S}}\delta\left(\mathbf{x}-\mathbf{x}^{(n)}\left(t\right)\right).
\end{equation}
 We note that the initial value problem consists of defining an initial
static distribution $P\left(\mathbf{x}\right)$ at $t=t_{0}$, and
then sampling it with initial points $\mathbf{x}^{(m)}\left(t_{0}\right)$.
For dynamics, the sampled paths are defined at $N_{T}$ discrete times
$t_{i}$, with a spacing $\Delta t$, such that exact ensemble average
values are recovered in the double limit of infinite sample and zero
step-size:

\begin{equation}
\left\langle o_{m}\left(t_{i}\right)\right\rangle =\lim_{\Delta t\rightarrow0}\lim_{N_{S}\rightarrow\infty}\frac{1}{N_{S}}\sum_{n=1}^{N_{S}}o_{m}\left(\mathbf{x}^{(n)}\left(t_{i}\right)\right).\label{eq:moments}
\end{equation}

These paths are obtained by numerical algorithms which propagate $\mathbf{x}^{(n)}(t_{i})$
to $\mathbf{x}^{(n)}(t_{i+1})$, with a truncation error vanishing
as $\Delta t\rightarrow0$. If we call $N_{S}$ the sample number,
the phase-space is made up of a set of distinct samples or trajectories.
These form an extended vector of parallel trajectories,
\begin{equation}
\mathbf{X}=\left(\begin{array}{ccc}
\mathbf{x}^{(1)T} & \cdots & \mathbf{x}^{(N_{S})T}\end{array}\right)^{T},
\end{equation}
of size $N_{S}d$. It is convenient to define the sampled observables
as a function $\bar{\mathbf{o}}$ of the extended vector of sampled
trajectories,
\begin{equation}
\bar{o}_{m}\left(\mathbf{X}\right)=\frac{1}{N_{S}}\sum_{n=1}^{N_{S}}o_{m}\left(\mathbf{x}^{(n)}\right).
\end{equation}

Throughout this paper, we will use the standard notation that all
extended vectors or matrices involving a dimension of the order of
the sample size $N_{S}$ are written in upper case like $\mathbf{X}$,
while objects that do not depend on the sample size are written in
lower case.

\subsection{Fokker-Planck and stochastic equations}

To explain the problem of interest here, we give a detailed analysis
of a stochastic differential equation or SDE\@. The underlying distribution
function $P\left(t,\mathbf{x}\right)$ satisfies the Fokker-Planck
equation~\cite{Gardiner1997},
\begin{equation}
\frac{\partial P\left(t,\mathbf{x}\right)}{\partial t}=\left[-\sum_{i}\frac{\partial}{\partial x_{i}}a_{i}(\mathbf{x})+\frac{1}{2}\sum_{ij}\frac{\partial^{2}}{\partial x_{i}\partial x_{j}}\mathsf{d}_{ij}(\mathbf{x})\right]P\left(t,\mathbf{x}\right)\,.\label{eq:FPE}
\end{equation}
We wish to minimize the errors for solving the corresponding It\={o}
SDE~\cite{Arnold1992-stochastic,Gardiner1997}, given in a standard
form by:
\begin{equation}
\mathrm{d}\mathbf{x}=\mathbf{a}\left(\mathbf{x}\right)\mathrm{d}t+\bm{\mathsf{b}}\left(\mathbf{x}\right)\mathrm{d}\mathbf{w},\label{eq:SDE}
\end{equation}
where $\bm{\mathsf{d}}=\bm{\mathsf{b}}\bm{\mathsf{b}}^{T}$, and the
Gaussian distributed real noise $\mathrm{d}\mathbf{w}$ has correlations
given by:
\begin{eqnarray}
\left\langle \mathrm{d}w_{i}\mathrm{d}w_{j}\right\rangle _{\infty} & = & \delta_{ij}\mathrm{d}t\nonumber \\
\left\langle \mathrm{d}w_{i}\right\rangle _{\infty} & = & 0.\label{eq:noise-averages}
\end{eqnarray}

Here the $\left\langle \dots\right\rangle $ notation means ensemble
average, so that
\begin{equation}
\left\langle \mathbf{x}\right\rangle _{N_{S}}\equiv\frac{1}{N_{S}}\sum_{n=1}^{N_{S}}x^{(n)},
\end{equation}
with $\langle\dots\rangle_{\infty}$ being the infinite ensemble limit.

\subsection{Sampling and error criteria}

To evaluate the computational error, we must define an error criterion
relevant to the entire sample ensemble $\mathbf{X}$. In practical
terms, this means one must calculate errors in each required observable,
rather than just the errors in a single trajectory. By minimizing
the total error over a finite number of samples, one should be able
to arrive at the most efficient algorithm that utilizes a given computational
resource.

To quantify the error we will use a weighted norm $\left\Vert P\right\Vert _{W}=\sqrt{\left(P,P\right)}\geq0$.
The integration error of a sampled, calculated stochastic distribution
$P_{S}$ relative to the exact distribution $P_{E}$ is then defined
as
\begin{equation}
\epsilon=\left\Vert P_{S}-P_{E}\right\Vert _{W},
\end{equation}
which depends on the choice of weighted norm.

The computed averages, $\left\langle o_{m}\right\rangle _{P_{S}}$,
are obtained through sampling, which means that typically one is interested
in calculating specific observables with minimum error. To evaluate
the computational accuracy relevant to these observables of interest,
we choose a particular set of observable quantities, $o_{m}$, for
$m=1,\ldots M$, to define the distribution norm, so that:
\begin{equation}
\left\Vert P\right\Vert _{W}^{2}\equiv\sum_{m=1}^{M}W_{m}\left|\left\langle o_{m}\right\rangle _{P}\right|^{2}.\label{eq:finite-moment-inner-product}
\end{equation}
Here $W_{m}$ is the relative weight assigned to observable $m$.
We assume the error is evaluated at a fixed time $t$, otherwise one
may wish to minimize the errors at each of a set of sample times.
In most of this paper, the observables are a finite set of one-dimensional
moments, $o_{m}\equiv x^{m}$, although other measures are certainly
possible. For simplicity, we choose $W_{m}=1$ from now on.

With this definition, the \emph{total} error is
\begin{equation}
\epsilon=\sqrt{\sum_{m=1}^{M}\left|\left\langle o_{m}\right\rangle _{P_{S}}-\left\langle o_{m}\right\rangle _{P_{E}}\right|^{2}},
\end{equation}

The integration error $\epsilon$ clearly must be calculated from
the entire vector of stochastic trajectories $\mathbf{x}$. The task
of defining an optimal stochastic integration method is to obtain
$\mathbf{x}\left(t\right)$ with a procedure that minimizes the total
error $\epsilon$ relative to some computational resource. Using the
definition that
\begin{equation}
\mu_{m}=\left\langle o_{m}\right\rangle _{P_{E}},
\end{equation}
we see that:

\begin{equation}
\epsilon=\sqrt{\sum_{m=1}^{M}\left|\bar{o}_{m}\left(\mathbf{X}\right)-\mu_{m}\right|^{2}}.\label{eq:distribution-error}
\end{equation}

Although different to the usual probability norms, we will use the
above definition throughout this paper, since it corresponds to the
operational requirements of a norm for a sampled distribution. For
a single observable, say $x$, $\epsilon^{2}$ is simply the usual
sample variance in $x$~\textemdash{} provided there are no other
errors present as well. However, there are usually other errors, including
errors due to arithmetic roundoff and the finite step-size of the
integration algorithm.

\subsection{Optimizing the sample number}

What is the optimal strategy to solve this SDE, given fixed total
computational resources and a goal of minimizing the total error?
With independent paths, the computational resource utilized is proportional
to $N=N_{S}N_{T}$, where $N_{S}$ is the number of samples and $N_{T}$
is the number of time-steps. One can invest total processing time
either in reducing the time-step, or in increasing the number of samples,
but reducing the step-size means using less samples, and vice-versa.

To understand these issues quantitatively, we first consider the optimal
strategy for error-reduction with traditional methods that involve
using independent sample paths. We will show that it is only important
to reduce the step-size or discretization error to the point where
it is some fraction of the sampling error.

Consider the total error for $N_{T}$ time-steps with a stochastic
integration method of global order $p$, that utilizes $N_{S}$ \emph{independent}
sample paths. From the central limit theorem, the global sampling
error $\epsilon_{S}$ due to the use of a finite sample, is
\begin{equation}
\epsilon_{S}=\sigma N_{S}^{-1/2},\label{eq:sampling-error}
\end{equation}
where $\sigma$ is the standard deviation of the measured quantity,
and the same type of scaling holds for sums over errors in multiple
observables.

There are additional errors $\epsilon_{T}$ over a fixed time interval
$T=N_{T}\Delta t$, due to the finite time-step $\mathrm{\Delta}t$.
If the local error of one time-step is $\epsilon_{\mathrm{\Delta}t}$,
then at best $\epsilon_{T}\approx N_{T}\epsilon_{\mathrm{\Delta}t}$.
The discretization order $p$ is therefore defined so that the global
time-step truncation error $\epsilon_{T}$ scales as
\begin{equation}
\epsilon_{T}=cN_{T}^{-p}\propto\Delta t^{p},\label{eq:discretization-error}
\end{equation}
where $c$ is an algorithm-dependent constant that depends on the
set of computed observables.

Since the two error sources are independent, the overall error is
given approximately by a sum of two terms:
\begin{equation}
\epsilon=\epsilon_{T}+\epsilon_{S}=cN_{T}{}^{-p}+\sigma N_{S}^{-1/2}.\label{eq:total-error}
\end{equation}
We wish to constrain the total processor time~\textemdash{} parallel
or otherwise~\textemdash{} so that $N_{S}N_{T}=N$ is bounded. The
total CPU time is then $T_{\mathrm{CPU}}=N\tau$, if one step in time
takes a real time duration $\tau$, and other overheads are negligible.

Next, consider how to optimize the tradeoff between the truncation
and sampling errors. Letting $\tilde{c}=cN^{-p}$, and keeping $N$
fixed, one obtains:
\begin{equation}
\epsilon=\tilde{c}N_{S}^{p}+\sigma N_{S}^{-1/2}
\end{equation}
This is minimized by choosing a combination of time-step and sample
number such that:

\begin{eqnarray}
N_{S} & = & N^{\frac{2p}{2p+1}}\left[\frac{\sigma}{2pc}\right]^{\frac{2}{2p+1}}\nonumber \\
N_{T} & = & N^{\frac{1}{2p+1}}\left[\frac{\sigma}{2pc}\right]^{-\frac{2}{2p+1}}
\end{eqnarray}

The total error is then:
\begin{equation}
\epsilon=\left[\frac{c\sigma^{2p}}{2pN^{p}}\right]^{\frac{1}{2p+1}}\left[1+\frac{1}{2p}\right]\label{eq:optimum-error}
\end{equation}

\subsection{Effective integration order}

From the analysis above, the best that one can do to reduce errors
is to obtain a scaling of $\epsilon\propto N^{-p_{\mathrm{eff}}}$,
where the effective order $p_{\mathrm{eff}}$ is
\begin{equation}
p_{\mathrm{eff}}=\frac{p}{2p+1}.
\end{equation}
This is because a higher order integration technique has no effect
on the sampling error, which eventually dominates the error calculation.
The best effective order, even with $p=\infty$, is $p_{\mathrm{eff}}=1/2$.
Given that stochastic higher order algorithms are complex and slow,
higher $p$ values are not always an advantage.

This optimal effective order requires the use of an optimal ratio
of samples to steps, which is:
\begin{equation}
\frac{N_{S}}{N_{T}}=N^{\frac{2p-1}{2p+1}}\left[\frac{\sigma}{2pc}\right]^{\frac{4}{2p+1}}\,.
\end{equation}
This ratio is difficult to calculate a-priori, since these constants
are not usually known in advance. One can estimate the sampling error
numerically by measuring the sample standard deviation $\sigma$,
which is known from the statistics of the simulation. The truncation
error constant $c$ is also measurable by changing the step-size.
Thus, both $c$ and $\sigma$ can be measured in software. An optimum
point is reached when the ratio between discretization and sampling
error is given by the extremely simple result that:
\begin{equation}
\frac{\epsilon_{T}}{\epsilon_{S}}=\frac{1}{2p}\,.
\end{equation}

In summary, it is only productive to reduce the step-size to the point
where the discretization error of $1/2p$ of the sampling error. After
this point is reached, one is better off to use more samples. Another
way to describe this is that the least computational cost for error
$\epsilon$ is
\begin{equation}
N\propto\epsilon^{-1/p_{\mathrm{eff}}}=\epsilon^{-\left(2p+1\right)/p}.
\end{equation}
There are multi-scale methods~\cite{Giles2008-multilevel} that can
achieve a cost of $N\propto\epsilon^{-2}\left(\log\epsilon\right)^{2}$,
which is a useful further improvement.

In summary, from the traditional analyses of ODE integration, one
might expect that the error would be reduced extremely rapidly with
a higher order method. This is not the case with an SDE.\emph{ }Because
the independent sampling error varies slowly with resource, this becomes
the dominant error with increasing resources.

\emph{Our conclusion is that sampling error strongly limits the effectiveness
of independent sampling methods for SDEs.}

\section{Initial conditions for parallel optimized sampling\label{sec:static:initial-conditions}}

We now describe a family of variance reduction methods using parallel
optimized sampling (POS), which locally eliminate sampling errors
for a finite set of observables $\mathbf{o}$, thus improving computational
efficiency. As described above, for $N_{S}$ samples, there is a $N_{S}d$-dimensional
extended vector of stochastic variables, $\mathbf{X}=\left(\begin{array}{ccc}
\mathbf{x}^{\left(1\right)T} & \cdots & \mathbf{x}^{\left(N_{S}\right)T}\end{array}\right)^{T}$, which is optimized so that $\mathbf{X}\rightarrow\mathbf{X}^{(\mathrm{opt})}$
. Both initial samples and the resulting stochastic trajectories are
optimized in parallel, thus making use of increasingly parallel hardware
capabilities.

Before deriving the POS algorithm, we need to address an issue arising
from the known initial condition for the SDE, with $\mathbf{x}$ having
a distribution $P\left(\mathbf{x}\right)$ at $t=t_{0}$. When we
derive the optimization conditions, one of the assumptions will be
that the moments of $\mathbf{x}$ are already optimized at the \emph{start}
of each time step, that is $\langle o_{m}\left(t\right)\rangle_{N_{S}}=\mu_{m}\left(t\right)$
for $m=1\dots M$. At every step of the SDE integration, if we obtain
the optimized stochastic vector from the previous step, then, since
we optimize the differentials of the moments, the resulting stochastic
variables will be optimized as well as possible, apart from any time-step
errors, which is a separate issue.

But how do we ensure the \emph{initial} set of moments of the distribution
$P\left(\mathbf{x}\right)$ is optimized at the first integration
step?

\subsection{Static POS}

We emphasize that our main goal is to minimize the variance for the
dynamical stochastic equations, treated next. As a consequence, we
do not use methods like Gaussian quadratures for the static variance
reduction problem. These are limited in the diversity of initial distributions
treated, and they cannot be easily generalized to the dynamic case.
Instead we take a set of random samples from $P\left(\mathbf{x}\right)$,
and optimize these initial sampled variables to reduce the variance
in moments of interest.

We call this static POS, as it does not involve dynamical time-evolution.
In other words, we wish to sample the initial values for the trajectories
\textbf{so that $M$ chosen observables are equal to their known exact
values}, that is,
\begin{equation}
\bar{\mathbf{o}}\left(\mathbf{X}\right)=\bm{\mu}\equiv\langle\mathbf{o}\rangle_{\infty},\label{eq:static-target-eqn}
\end{equation}
where $\boldsymbol{\mu}$ is an $M$-dimensional vector. That is,
we wish to set $\epsilon_{P}=\epsilon_{P}(t_{0})=0$, where $\epsilon_{P}$
is the total error~(\ref{eq:distribution-error}) in the initial
distribution, using an observable-based error measure. In the examples
given in this paper, we treat all moments up to a given order. More
generally, at high dimension $d$, treating all correlations or moments
to a fixed order is exponentially hard. Therefore, it is necessary
to select moments according some importance criterion, the effects
of which will be treated in a later publication.

This is a classic example of nonlinear, multidimensional root-finding,
with optimization over the dimensionality of the full space of samples.
It can be a challenge to solve these equations efficiently, with resources
no greater than $\mathcal{O}(N_{S})$, and with minimal changes to
the original sampled estimate for $\mathbf{X}(t_{0})$, which is labelled
$\mathbf{X}_{(0)}$. However, in examples treated here, we find that
the numerical problem can be very accurately solved, naturally with
machine-limited numerical precision. Here we use an iterative, modified
Newton-Raphson method for the solution. Other techniques are available
as well.

In an iterative Newton-Raphson approach, we use stochastic methods
for the first estimate, $\mathbf{X}_{(0)}$. At each iteration, we
set $\mathbf{X}\rightarrow\mathbf{X}+\delta\mathbf{X}$ to obtain
the next estimate $\mathbf{X}_{(i+1)}$, which has a better fit to
the required moments, so that $\mathbf{X}_{(i+1)}=\mathbf{X}_{(i)}+\delta\mathbf{X}_{(i+1)}$,
where $\|\delta\mathbf{X}_{(i+1)}\|\ll\|\boldsymbol{X}_{(i)}\|$,
and $\|\dots\|$ denotes the Euclidean norm.

Expanding to first order in $\delta\mathbf{X}_{(i+1)}$, the conditions
we require are:

\begin{eqnarray}
\bar{\mathbf{o}}^{\mathrm{(opt)}} & = & \bm{\mu}\nonumber \\
 & = & \bar{\mathbf{o}}\left(\mathbf{X}_{(i)}\right)+\bm{\mathsf{J}}\left(\mathbf{X}_{(i)}\right)\delta\mathbf{X}_{(i+1)}+\mathcal{O}\left(\delta\mathbf{X}_{(i+1)}^{2}\right),
\end{eqnarray}
where $\bm{\mathsf{J}}$ is a Jacobian of $\bar{\mathbf{o}}$:
\begin{equation}
\mathsf{J}_{mn}\left(\mathbf{X}\right)=\frac{\partial\bar{o}_{m}\left(\mathbf{X}\right)}{\partial X_{n}}\,.\label{eq:Derivative matrix}
\end{equation}

In matrix form, the iteration equation is therefore

\begin{equation}
\bm{\mathsf{J}}_{(i)}\delta\mathbf{X}_{(i+1)}=\bm{\mu}-\bar{\mathbf{o}}_{(i)}.\label{eq:iteration-equation}
\end{equation}
where we use the obvious notation that $\bar{\mathbf{o}}_{(i)}\equiv\bar{\mathbf{o}}\left(\mathbf{X}_{(i)}\right)$
and $\bm{\mathsf{J}}_{(i)}\equiv\bm{\mathsf{J}}\left(\mathbf{X}_{(i)}\right)$.
We note that out method is scale-invariant in the sense that it will
result in the exact same solutions if we rescale any observable and
its exact value as $o_{i}(\mathbf{X})\rightarrow c_{i}o_{i}(\mathbf{X})$,
where $c_{i}$ is a constant.

Each step of the iterative solution is underdetermined, since the
relevant matrices are not square, and hence not invertible. Such problems
are relatively common in optimization, and here we analyze two different
techniques to solve them.

\subsection{One-dimensional example of iteration equations}

We will use as an illustrative example the situation when the observables
are moments of a one-dimensional stochastic equation. In this example
$o_{m}=x^{m}$, so the conditions are:

\begin{eqnarray}
\left[\mathbf{X}^{\mathrm{(opt)}}\right]^{m}\cdot\mathbf{1} & = & N_{S}\mu_{m}\nonumber \\
 & = & \mathbf{X}_{(i)}^{m}\cdot\mathbf{1}+m\mathbf{X}_{(i)}^{m-1}\cdot\delta\mathbf{X}_{(i+1)}+\mathcal{O}\left(\delta\mathbf{X}_{(i+1)}^{2}\right),
\end{eqnarray}
where the powers of matrices and vectors are understood in the element-wise
(Hadamard) sense (e.g., $\mathbf{X}^{2}\equiv\mathbf{X}\circ\mathbf{X}$).
Hence, in this case the observable vector and matrix of derivatives
are:

\begin{eqnarray}
\bar{\mathbf{o}}\left(\mathbf{X}\right) & = & \frac{1}{N_{S}}\sum_{n=1}^{N_{S}}\left(\begin{array}{c}
x^{(n)}\\
\cdots\\{}
[x^{(n)}]^{M}
\end{array}\right),\nonumber \\
\bm{\mathsf{J}}\left(\mathbf{X}\right) & = & \frac{1}{N_{S}}\left(\begin{array}{ccc}
1 & \cdots & 1\\
 & \cdots\\
M\left[x^{(1)}\right]^{M-1} & \cdots & M\left[x^{(N_{S})}\right]^{M-1}
\end{array}\right).\label{eq:XY-matrices}
\end{eqnarray}

\subsection{Static iteration solution with a matrix pseudo-inverse}

The simplest method to solve the iteration equation~(\ref{eq:iteration-equation})
is to use a matrix pseudo-inverse. This has another advantage, since
the pseudo-inverse has the property that it generates a least-squares
solution, with minimal changes, as required for this algorithm.

Thus, we have the iterative procedure

\begin{equation}
\mathbf{X}_{(i+1)}=\mathbf{X}_{(i)}+\bm{\mathsf{J}}_{(i)}^{+}\left(\bm{\mu}-\bar{\mathbf{o}}_{(i)}\right).
\end{equation}
where $\bm{\mathsf{J}}^{+}$ indicates the pseudo-inverse of the $M\times N_{S}$
Jacobian $\bm{\mathsf{J}}$.

This pseudo-inverse method is simple, general and effective. It satisfies
our requirements given above. However, standard pseudo-inverse software
using singular value decomposition (SVD) methods has drawbacks. It
is relatively slow, can have stability problems and may not be available
on all programming platforms.

\subsection{Improved matrix inversion iterations\label{sub:static-pos:expansion}}

To improve efficiency and overcome problems with the use of general
purpose pseudo-inverses, there is a faster alternative method for
solving~(\ref{eq:iteration-equation}).

The pseudo-inverse is a limit. If $\bm{\mathsf{J}}^{\dagger}$ denotes
the conjugate transpose of the matrix $\bm{\mathsf{J}}$, then:
\begin{equation}
\bm{\mathsf{J}}^{+}\equiv\lim_{\alpha\rightarrow0}\bm{\mathsf{J}}^{\dagger}\left(\bm{\mathsf{J}}\bm{\mathsf{J}}^{\dagger}+\alpha\mathbf{I}\right)^{-1}.
\end{equation}
Hence, provided $\bm{\mathsf{u}}=\bm{\mathsf{J}}\bm{\mathsf{J}}{}^{\dagger}$
is invertible, one can directly obtain that $\bm{\mathsf{J}}^{+}=\bm{\mathsf{J}}^{\dagger}\bm{\mathsf{u}}^{-1}$.
It is important to note here that $\bm{\mathsf{u}}$ is only an $M\times M$
matrix, regardless of how many samples were used originally.

This gives as many unknown coefficients as there are equations. As
a result, we only need solve $M$ linear equations in $M$ unknowns,
for which there are efficient techniques, resulting in the iteration
procedure

\begin{equation}
\mathbf{X}_{(i+1)}=\mathbf{X}_{(i)}+\bm{\mathsf{J}}_{(i)}^{\dagger}\bm{\mathsf{u}}_{(i)}^{-1}\left(\bm{\mu}-\bar{\mathbf{o}}_{(i)}\right).\label{eq:static iteration algorithm}
\end{equation}

Here the matrix $\bm{\mathsf{u}}$ being inverted is an $M\times M$
positive-semidefinite matrix. We use the lower case for it to indicate
that it does not take an extended dimensionality. It is generally
invertible and relatively small compared to the sample size, resulting
in a highly efficient optimization algorithm. Numerical investigations
show that this simple method is convergent for large sample number
$N_{S}$, which corresponds to cases where the initial sample has
moments close to their ideal values.

As the stopping condition in the numerical examples given later, we
choose the condition that the changes are small relative to the modulus
of the current $\mathbf{X}$ value, i.e. $\|\delta\mathbf{X}_{(i+1)}\|/\|\mathbf{X}_{(i)}\|<\eta=10^{-8}$,
where $\eta^{2}$ is the numerical machine precision, as previously.
An iteration limit of $I_{\max}=50$ was used, although as we see
later, the typical iteration number is much lower than this.

We note that this simplified method is only useful when $\bm{\mathsf{u}}$
is invertible, otherwise the full SVD techniques or other more robust
algorithms should be used. This is not a major limitation: as long
as the first derivatives of our observables are linearly independent,
$\bm{\mathsf{u}}$ will be invertible. As often the case with such
numerical methods, performance can be further improved through linear
scaling or translation operations.

However, sometimes the Newton-Raphson method doesn't converge, not
because $\bm{\mathsf{u}}$ lacks invertibility (which is rare), but
because the initial guess is too far away from the solution. This
can be detected by checking if the norm of $\bm{\Delta}$ increases.
In such cases, one can back-track and try again with a different initial
estimate. In the case of time-evolution dynamics, treated next, a
smaller time-step can also be used to ensure the initial guess is
close to hte solution.

\subsection{One-dimensional case: power series method}

As an equivalent way to understand our approach in the one-dimensional
case, we can expand $\delta\mathbf{X}$, which is the iterative \textbf{change}
in $\mathbf{X}$, using the derivative matrix and a vector of changes
in observables, $\delta\mathbf{o}$. With observables as moments,
this is a series in $\mathbf{X}^{\dagger}$ up to order $M-1$, with
$M$ coefficients written as a vector $\bm{\delta}$. The expansion
is given by:

\begin{equation}
\delta\mathbf{X}=\bm{\mathsf{J}}^{\dagger}\delta\bar{\mathbf{o}}=\frac{1}{N_{S}}\sum_{m=1}^{M}m\left(\mathbf{X}^{\dagger}\right)^{m-1}\delta\bar{o}_{m}.
\end{equation}

Defining a square, positive definite matrix:
\begin{equation}
\bm{\mathsf{u}}_{(i)}=\bm{\mathsf{J}}_{(i)}\bm{\mathsf{J}}_{(i)}^{\dagger},
\end{equation}
the basic iteration equation~(\ref{eq:iteration-equation}) becomes:
\begin{equation}
\bm{\mathsf{u}}_{(i)}\delta\bar{\mathbf{o}}{}_{(i+1)}=\bm{\mu}-\bar{\mathbf{o}}_{(i)}.
\end{equation}

Since the solution for $\bm{\delta}{}_{(i+1)}$ requires the inversion
of $\bm{\mathsf{u}}_{(i)}$, this \emph{leads to the same final algorithm}
that is defined above, in Eq~(\ref{eq:static iteration algorithm}).

\section{Numerical examples: Static optimization\label{sec:parallel-optimized-sampling}}

In this section, we perform numerical optimizations of initial conditions,
following the methods of the previous section. This is important because
the initial sampling error can dominate the sampling error found throughout
the calculation, if it is not optimized. For definiteness, in most
of the numerical examples of POS algorithms analyzed here, the $m$-th
observable is a finite moment of a one-dimensional real distribution.
Hence,
\begin{equation}
\bar{o}_{m}\equiv\left\langle o_{m}\right\rangle _{N_{S}}=\frac{1}{N_{S}}\sum_{n=1}^{N_{S}}\left[x^{(n)}\right]^{m}=\frac{\mathbf{X}^{m}\cdot\mathbf{1}}{N_{S}},
\end{equation}
where $\mathbf{1}$ is an $N_{S}$-dimensional unit column vector.
We optimize moments, for $m=1,\dots M$. As well as calculating the
resulting error in these optimized moments, we also obtain the errors
in other moments. An important issue is that the non-optimized moments
should not have increased errors. We will also consider the ensuing
results for two less regular observables: $o_{\exp}(x)=\exp(x)$ (a
combination of an infinite number of moments) and $o_{||}(x)=|x|$
(an irregular function).

\subsection{Static optimization\label{sub:benchmarks:static}}

In this section we will test the performance of the iterative algorithm
from Section~\ref{sec:static:initial-conditions} applied to the
most common initial condition~\textemdash{} a normal distribution.
However, we emphasize that this choice is purely for convenience,
and any initial distribution could be used. In this section we will
use $\mu=\mu_{1}=0$ for all the tests and set the variance $\sigma^{2}=\mu_{2}=1$.
This helps to avoid precision loss when calculating central moments,
which is important, since the method is extremely efficient and optimizes
the moments almost to the limit of numerical precision, as seen below.
If maximum precision is required for cases when $\mu\ne0$, one can
subtract the mean of the initial distribution so that $\mu=0$ in
the new variables.

With $\mu=0$, the moments of the test distribution are given by

\begin{equation}
\mu_{m}=\left\{ \begin{array}{cc}
0, & m\ \mathrm{is\ odd},\\
\sigma^{m}\left(m-1\right)!!, & m\ \mathrm{is\ even}.
\end{array}\right.
\end{equation}
The number of moments being optimized here is $M=6$, but this is
not essential. For the additional observables the expected values
are the expectation of the log-normal distribution

\begin{equation}
\mu_{\exp}=\exp\left(\frac{\sigma^{2}}{2}\right),
\end{equation}
and the expectation of the folded normal distribution

\begin{equation}
\mu_{||}=\sigma\sqrt{\frac{2}{\pi}}.
\end{equation}

We test the behavior of the algorithms in $10000$ independent optimization
attempts, each with $N_{S}=1000$ parallel samples, in order to collect
statistics. In each attempt, a $N_{S}$-dimensional vector $\mathbf{X}$
is generated using random sampling and the matrix inversion iterative
method from Section~\ref{sub:static-pos:expansion} is applied to
it. We then plot the distribution of measures of interest, which indicate
the performance of the method. This demonstrates how well the algorithm
performs when tested over a range of different initial random samples,
each drawn from the same underlying normal distribution.

\subsection{Distribution of the optimized variance reduction}

There are several quantities we are interested in. We first investigate
the distribution of the optimized variance. This gives an idea of
not just how much the variance is reduced, but also the probability
of a given amount of variance reduction.

\begin{figure}
\includegraphics{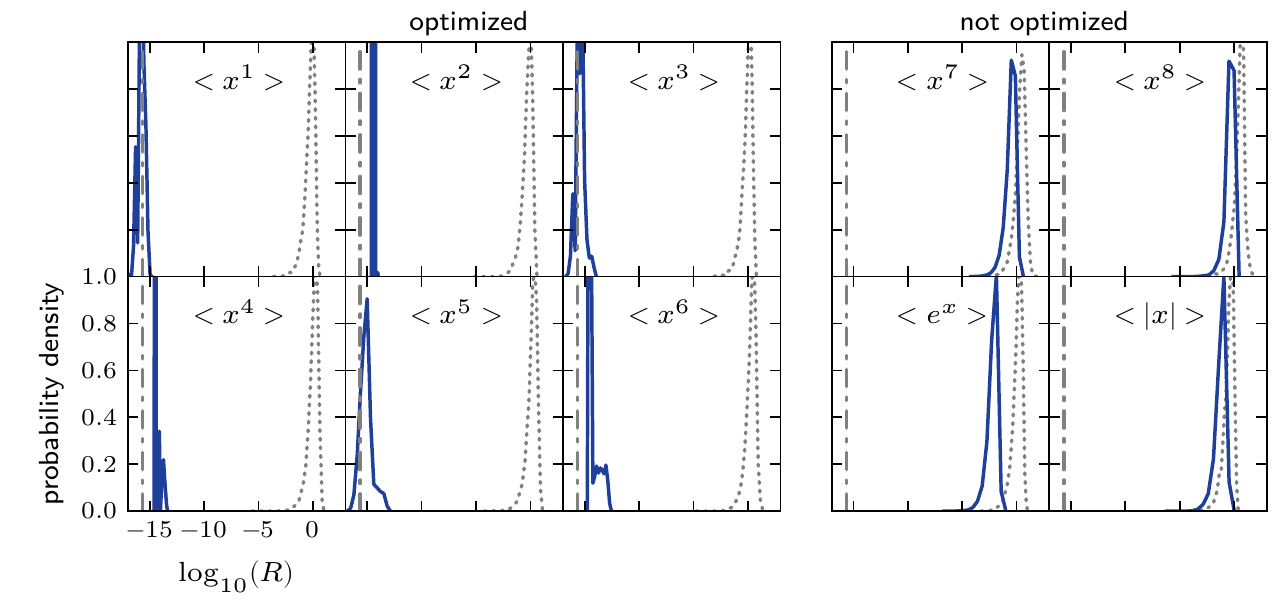}

\protect\caption{\label{fig:static:errors}Distribution of normalized errors in the
static optimization tests for the matrix inversion method with $\sigma=1$,
$N_{S}=1000$ (solid blue lines). The grey dotted curves show the
normalized errors of a randomly sampled vector. The dash-dotted line
denotes the numerical precision limit ($\epsilon_{\mathrm{num}}\approx2.22\times10^{-16}$
for double precision floating point numbers). Zero error results are
removed from binning for ease of plotting.}
\end{figure}

In more detail, we plot the numerically obtained distribution of final
central moments of the vector compared to the desired values:

\begin{equation}
\tilde{R}_{m}\left(\mathbf{X}\right)=\left|\frac{1}{N_{S}}o(\mathbf{X})\cdot\mathbf{1}-\mu_{m}\right|,
\end{equation}
where we will compare this quantity for initial randomly sampled vectors
and optimized vectors. It is convenient to scale this to the expectation
of the error in the sampled moment (for random samples) for the normal
distribution:

\begin{equation}
R_{m}=\frac{\tilde{R}_{m}}{\sigma^{m}\sqrt{m!/N_{S}}},
\end{equation}
so the average of $R_{m}$ for the initial randomly sampled vectors
will be close to $1$ for all $m$. Fig.~\ref{fig:static:errors}
shows this quantity both for the moments we are optimizing ($m=1\dots M$)
and for two higher moments that are not optimized. The errors for
the additional observables are normalized simply as

\begin{equation}
R_{\exp,||}=\frac{\tilde{R}_{\exp,||}}{\sqrt{N_{S}}}.
\end{equation}

The optimization procedure was able to reduce errors almost to the
limit of numerical precision of $R_{m}\sim10^{-15}$ in many cases.
The moments that were not directly optimized show some improvement
as well, but typically closer to a single order of magnitude. Since
we use logarithmic graphs, and $R_{m}$ can sometimes be zero, we
excluded these zero results from the binning.

\subsection{Number of iterations and optimization distance}

It is also useful to know the number of iterations performed in each
of the algorithms. This gives an idea of the numerical cost, which
we investigate in greater detail below. If too many iterations are
needed, the time taken will increase to the point that one may as
well use more traditional methods.

\begin{figure}
\includegraphics{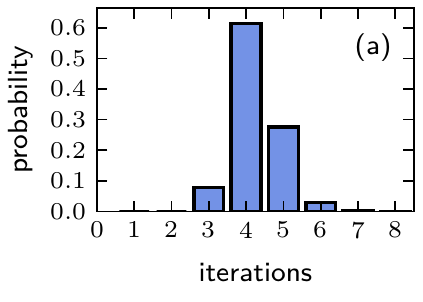}\includegraphics{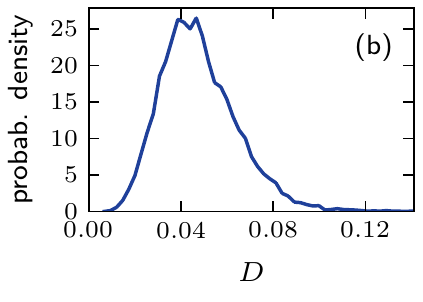}\includegraphics{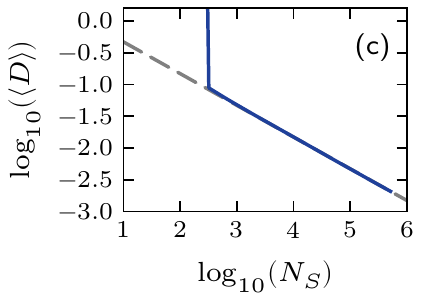}

\protect\caption{\label{fig:static:distance-and-iterations}(a) Number of iterations
taken by the matrix inversion method, $\sigma=1$, $N_{S}=1000$.
(b) Relative distance $D$ between the optimized vector and the initial
randomly sampled vector in the static optimization tests for the matrix
inversion method, $\sigma=1$, $N_{S}=1000$. (c) Average relative
distance $\langle D\rangle$ depending on the number of trajectories
$N_{S}$ in the static optimization tests for the matrix inversion
method, $\sigma=1$. The grey dashed line denotes the slope of a $1/\sqrt{N_{S}}$
dependence. }
\end{figure}

Fig.~\ref{fig:static:distance-and-iterations}(a) demonstrates that
only a small number of iterations~\textemdash{} typically four~\textemdash{}
are necessary for the method to converge, given a reasonably large
ensemble. The third quantity of interest is the relative distance
of the optimized vector from the initial randomly sampled one:

\begin{equation}
D=\frac{\|\mathbf{X}^{(\mathrm{opt})}-\mathbf{X}\|}{\|\mathbf{X}\|}.
\end{equation}
We would like this to be much less than $1$, which means that the
optimized vector is still approximately ``random''. Fig.~\ref{fig:static:distance-and-iterations}(b)
shows that it is indeed the case, and the average deviation is $D\approx5\times10^{-3}$
for $\sigma=1$. Fig.~\ref{fig:static:distance-and-iterations}(c)
shows that the results have the desirable property that the relative
distance between the optimized and non-optimized trajectories scales
with $1/\sqrt{N_{S}}$. This is expected from the fact that the initial
sampled moment relative errors scale as $1/\sqrt{N_{S}}$, hence one
expects the minimum required change to be of this order.

However, if there are too few trajectories, the simple matrix inversion
iteration method diverges, producing unacceptably large changes $D$
in the distribution. For these small $N_{S}$ values, the iteration
limit of $I_{\max}=50$ was reached, indicating a lack of convergence.

\begin{figure}
\includegraphics{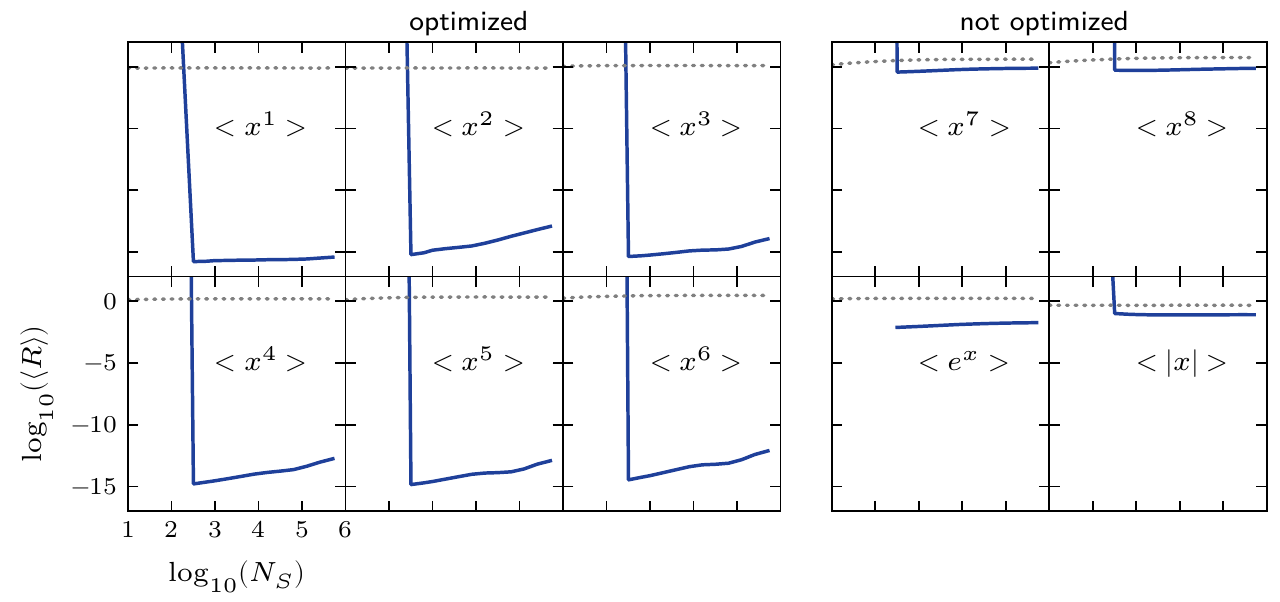}

\protect\caption{\label{fig:static:error-vs-N}Dependence of normalized errors on the
number of trajectories $N_{S}$ in the static optimization tests for
the expansion method, $\sigma=1$. The dotted curve shows the dependence
for randomly sampled vectors.}
\end{figure}

To investigate this more carefully, the dependence of the final error
$R$ on the number of trajectories $N_{S}$, is plotted in Fig.~\ref{fig:static:error-vs-N}.
The results also display a sharp cut-off in the number of trajectories
below which the iterative approach does not converge. More robust
minimization algorithms could presumably solve this, but the issue
was not investigated here for space reasons.

\subsection{Initial optimization cost}

For practical applications, the time taken for the optimization is
also important. More accurately, the quantity of interest is the efficiency~\cite{L'Ecuyer1994-efficiency},
which in our case can be expressed as $\mathrm{Eff=1/\left(T_{\mathrm{CPU}}\tilde{R}^{2}\right)}$,
where $T_{\mathrm{CPU}}$ is the time required to generate the vector
$\mathbf{X}^{(\mathrm{opt})}$.

\begin{figure}
\includegraphics{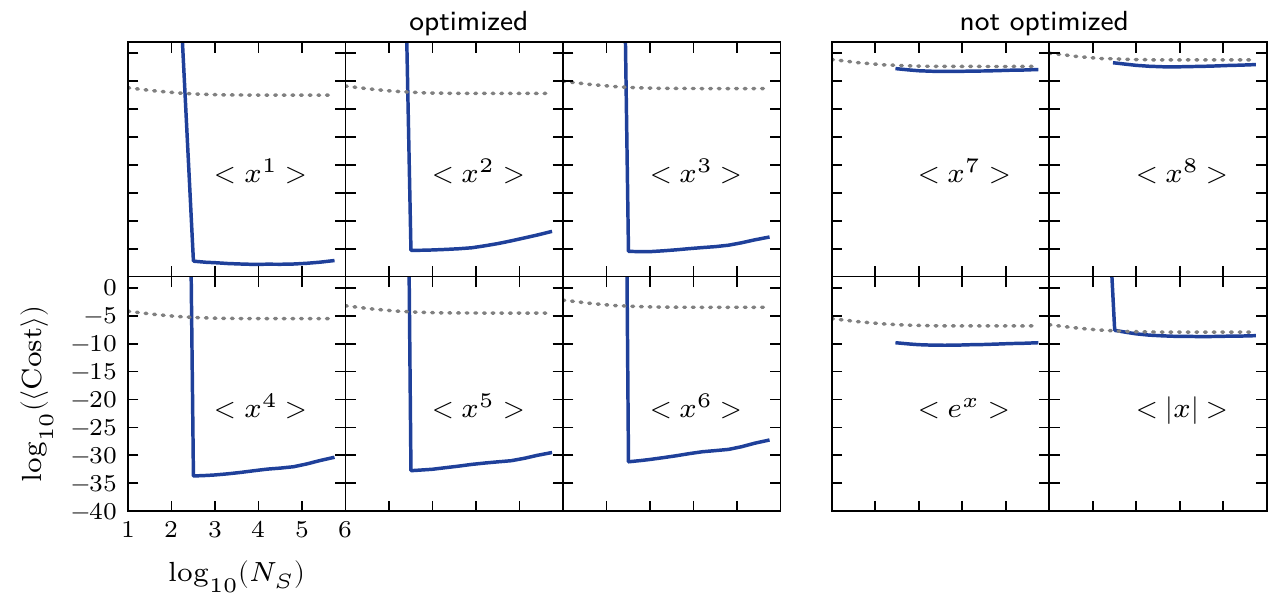}

\protect\caption{\label{fig:static:efficiency}Dependence of the cost on the number
of trajectories $N_{S}$ in the static optimization tests for the
expansion method, $\sigma=1$. The dotted curve shows the dependence
for randomly sampled vectors.}
\end{figure}
Since in our case in a significant fraction of optimization attempts
$\tilde{R}$ for optimized moments is exactly zero, we plot the cost

\begin{equation}
\mathrm{Cost}=\frac{1}{\mathrm{Eff}}\equiv T_{\mathrm{CPU}}\tilde{R}^{2}.
\end{equation}
averaged over attempts. Fig.~\ref{fig:static:efficiency} shows this
quantity for all investigated observables. For the optimized moments
($m=1\dots6$) the cost of POS is significantly lower than the cost
of the purely random sampling (starting from a certain critical number
of trajectories). The situation is different for non-optimized observables,
where POS still demonstrates a decreased cost, but the difference
is not as overwhelming.

In conclusion, the expansion iterative method performs well for a
wide range of parameters, while being fast and easy to implement.
Convergence is extremely rapid when the number of stochastic trajectories
is increased above a critical value. However, we note that many other
methods of nonlinear root-finding exist, and we do not exclude the
possibility of better techniques being available.

We use this approach in the SDE examples given below, to pre-optimize
the starting conditions for the SDE integration tests.

\section{Parallel dynamic optimization\label{sec:dynamic}}

The general approach here is to assume that samples with an optimized
set of observables are known at a given time $t$. Given this, a noise
ensemble is chosen to give the lowest possible errors in the observables
at a later time $t+\Delta t$. Overall convergence is then achieved
by taking a large $N_{S}$ limit. The goal of the algorithm is that
the observable moments are obtained with minimal error. For efficiency,
the resulting calculation complexity should be linear in $N_{S}$.

In summary, the POS algorithm must therefore meet three essential
requirements:
\begin{enumerate}
\item Equations for a finite set of $M$ observables are optimized both
initially, and for local changes in time to a given order in $\Delta t$.
\item Non-optimized observables still have their correct values to the same
order in $\Delta t$ in the limit of $N_{S}\rightarrow\infty$.
\item The computational cost of the algorithm should be no worse than $\mathcal{O}(N_{S})$.
\end{enumerate}
If the observables chosen are moments, then, since higher order moments
depend on lower order ones, these will experience a degree of optimization
as well, even though not directly optimized. POS methods can also
be combined with higher-order time-step methods, but here we focus
on sampling error reduction for ease of explanation. Other types of
observable also experience a degree of variance reduction as well.
Typically this is not as great as experienced by the optimized moments,
but this depends on the observable and the stochastic equation.

There is an analogy between this approach and the method of moment
hierarchies in statistical physics~\cite{Risken1996}. The difference
is that rather than an arbitrary truncation of the moment hierarchy,
the higher order moments are estimated in an unbiased way via sampling.

We now describe how to extend this initial optimization algorithm
to treat dynamical optimization of moments during stochastic time
evolution of $\mathbf{x}$. This involves stochastic noise terms and
deterministic drift terms.

\subsection{Euler-Maruyama algorithm}

For simplicity, we only treat the optimization of the Euler-Maruyama
integration scheme for an SDE in It\={o} form~\cite{Kloeden1992},
which has truncation order $p=1$ for convergence of moments. This
is a discrete expression of the standard Ito SDE, so that for a finite
step in time $\Delta t$,

\begin{equation}
\Delta\mathbf{x}=\mathbf{a}\Delta t+\bm{\mathsf{b}}\Delta\mathbf{w}.
\end{equation}

Next, we introduce $\mathbf{A}=\left(\begin{array}{ccc}
\mathbf{a}^{T}\left(\mathbf{x}^{(1)}\right) & \cdots & \mathbf{a}^{T}\left(\mathbf{x}^{\left(N_{S}\right)}\right)\end{array}\right)^{T}$ as an extended vector of drift coefficients, $\bm{\mathsf{B}}=\mathrm{diag}\left(\begin{array}{ccc}
\bm{\mathsf{b}}\left(\mathbf{x}^{(1)}\right) & \cdots & \bm{\mathsf{b}}\left(\mathbf{x}^{\left(N_{S}\right)}\right)\end{array}\right)$ (a block-diagonal matrix with $\bm{\mathsf{b}}\left(\mathbf{x}^{(i)}\right)$
elements) as an extended matrix of noise terms, and $\Delta\mathbf{W}=\left(\begin{array}{ccc}
\left(\Delta\mathbf{w}^{(1)}\right)^{T} & \cdots & \left(\Delta\mathbf{w}^{\left(N_{S}\right)}\right)^{T}\end{array}\right)^{T}$ as an extended vector of noises. The parallel sample vector $\mathbf{X}$
will satisfy the following vector SDE, now of dimension $N_{S}d$:

\begin{equation}
\Delta\mathbf{X}=\mathbf{A}\Delta t+\bm{\mathsf{B}}\Delta\mathbf{W}.\label{eq:parallel-trajectory equation}
\end{equation}

One requires that the sampled observables are given as closely as
possible by the known infinite ensemble results, that is,

\begin{equation}
\bar{\mathbf{o}}\left(\mathbf{X}+\mathrm{\Delta}\mathbf{X}\right)-\bar{\mathbf{o}}\left(\mathbf{X}\right)=\left\langle \Delta\left(\mathbf{o}(\mathbf{x})\right)\right\rangle _{\infty}+\mathcal{O}(\Delta t^{2}),\label{eq:moment-equations}
\end{equation}
where $\left\langle \Delta\left(\mathbf{o}(\mathbf{x})\right)\right\rangle _{\infty}$
is the ideal observable change for an infinite ensemble to order $\Delta t$.
This is given either by an application of a single-step approximation
to the Fokker-Planck equation~(\ref{eq:FPE}) followed by partial
integration, or equivalently by an application of It\={o}'s formula~\cite{Gardiner1997}
to the stochastic equation~(\ref{eq:SDE}):

\begin{equation}
\left\langle \Delta o_{m}\left(\mathbf{x}\right)\right\rangle _{\infty}=\left\langle \frac{\partial o_{m}}{\partial x_{i}}a_{i}+\frac{1}{2}\frac{\partial^{2}o_{m}}{\partial x_{i}\partial x_{j}}\mathsf{d}_{ij}\right\rangle _{\infty}\Delta t+\mathcal{O}\left(\Delta t^{2}\right)\label{eq:ideal-observable-changes-appr-dt}
\end{equation}

In order to be able to calculate these moments in practice we must
make another approximation, since we only know the required averages
for a finite ensemble of trajectories. From the central limit theorem,
if the individual variances are finite then:

\begin{equation}
\left\langle \frac{\partial o_{m}}{\partial x_{i}}a_{i}+\frac{1}{2}\frac{\partial^{2}o_{m}}{\partial x_{i}\partial x_{j}}\mathsf{d}_{ij}\right\rangle _{\infty}=\left\langle \frac{\partial o_{m}}{\partial x_{i}}a_{i}+\frac{1}{2}\frac{\partial^{2}o_{m}}{\partial x_{i}\partial x_{j}}\mathsf{d}_{ij}\right\rangle _{N_{S}}+\mathcal{O}\left(1/\sqrt{N_{S}}\right).\label{eq:ideal-observable-changes-appr-N}
\end{equation}
In some cases the equality is exact; see the discussion of the error
introduced by this approximation in Section~\ref{sub:dynamic:error-propagation}.

\subsection{One-dimensional case}

To illustrate how the method works for a single variable, we note
that in the one-dimensional moment-based case, this result reduces
to:
\begin{eqnarray}
\left\langle \Delta\left(x^{m}\right)\right\rangle _{\infty} & \approx & m\left\langle x^{m-1}a+\frac{m-1}{2}x^{m-2}\mathsf{b}^{2}\right\rangle _{N_{S}}\Delta t.\label{eq:ideal moment changes}
\end{eqnarray}
These methods can also be extended to higher orders in the step-size,
but here we treat the equations up to the first order in $\Delta t$,
and truncate higher orders for simplicity. Higher order algorithms
will be treated elsewhere.

At this stage, we point out the existence of more than one possible
strategy for satisfying the moment equations. Since each has its own
distinct advantages and disadvantages, two particular approaches will
be treated here. They are described in the following two subsections.

Both of the resulting POS methods have the following useful features:
\begin{itemize}
\item Uniqueness~\textemdash{} the matrix iteration equations have unique
results.
\item Parameter independence~\textemdash{} only the observables that are
optimized are specified, not an arbitrary parameter.
\item Linear complexity~\textemdash{} the time taken scales as $N_{S}$,
since the inner products only require $N_{S}$ operations.
\end{itemize}

\subsection{Error propagation\label{sub:dynamic:error-propagation}}

With any dynamical POS algorithm, the goal of the optimization is
to remove any difference between the sampled observables and the infinite
trajectory observables. This can be achieved exactly for static optimization,
but there are error propagation effects to be considered in the dynamical
case. Before treating the step-wise optimization method, we briefly
consider the potential effect of error propagation.

To understand this, we note that error propagation effects depend
on both the choice of observables and the structure of the equations
themselves. Ideally, the optimized set of observables form a closed
set of equations, but this is rarely the case in practical nonlinear
calculations.

We illustrate the effects of error propagation by considering the
one-dimensional moment-based equations. Initially, the error in the
$p$-th moment is:

\begin{equation}
\epsilon_{p}=\left|\langle x^{p}+\mathrm{\Delta}x_{(0)}^{p}\rangle_{N_{S}}-\langle x^{p}+\mathrm{\Delta}x^{p}\rangle_{\infty}\right|\ne0,
\end{equation}
which varies as $1/\sqrt{N_{S}}$ for random noises. The aim of our
optimization is to set the difference to zero, apart from time-step
errors. Can this be achieved exactly, assuming that the stochastic
variables themselves are initially optimized at the start of the simulation,
that is, $\langle x^{p}\rangle_{N_{S}}=\langle x^{p}\rangle_{\infty}$
for $p=1,\dots M$?

To answer this, note that the drift and the diffusion terms can be
expanded in the powers of $x$ as

\begin{equation}
a(x,t)=\sum_{j=0}^{\infty}f_{j}(t)x^{j},\quad b^{2}(x,t)=\sum_{j=0}^{\infty}g_{j}(t)x^{j}.
\end{equation}
After substituting the expansions into the condition expressions,
and expanding to order $(\Delta t)$,

\begin{eqnarray}
\epsilon_{p}= &  & p^{2}\left|\left\langle x^{p-1}\Delta w+\frac{p-1}{2}x^{p-2}\left(\mathrm{\Delta}w^{2}-\mathsf{b}^{2}\mathrm{\Delta}t\right)\right\rangle _{N_{S}}\right.\nonumber \\
 &  & +\sum_{j=M-p+2}^{\infty}f_{j}(t)\left(\langle x^{p+j-1}\rangle_{N_{S}}-\langle x^{p+j-1}\rangle_{\infty}\right)\mathrm{\Delta}t\label{eq:total error-1}\\
 &  & \left.+\frac{p-1}{2}\sum_{j=M-p+3}^{\infty}g_{j}(t)\left(\langle x^{p+j-2}\rangle_{N_{S}}-\langle x^{p+j-2}\rangle_{\infty}\right)\Delta t\right|^{2}.\nonumber
\end{eqnarray}
The last two terms depend on unoptimized and thus unknown differences
$\langle x^{j}\rangle_{N_{S}}-\langle x^{j}\rangle_{\infty}$ for
$j>M$. In deriving the POS algorithm, we assume that all the moments
of the sampled distribution are correct at the start of the step in
time, which allows us to neglect these higher-order moment differences.

One possibility is that the coefficients $f$ and $g$ for the corresponding
indices are in fact zero. This is equivalent to having $a$ at most
linear in $x$, and $b^{2}$ at most quadratic in $x$. If this condition
is satisfied, then higher order moments remain uncoupled to lower
order moments, and no error propagation occurs apart from the usual
error propagation in an ODE method.

If this condition is not satisfied, neglecting these terms will result
in unoptimized high order moments ``leaking in'' to low order moments
over time, with the effect more noticeable for orders closer to $M$.
We ignore this effect in deriving the algorithm, which means that
we can only \emph{locally} remove sampling errors exactly, as we show
later. As a result, higher-order moment errors gradually degrade the
optimized moments. This results in globally increased errors for nonlinear
SDE solutions, which is demonstrated in the numerical examples described
in the next section. The consequence is that it is no longer possible
to reach machine precision. However, substantial variance reduction
is still possible.

\subsection{Combined~optimization}

In combined optimization, all time orders in the sampled moments are
combined together, to give a moment estimate in an analogous form
to that treated previously for the initial conditions. The advantage
of this approach is its ease of implementation and simplicity.

However, unlike the case of the initial distribution, we typically
do not know the infinite ensemble average moments exactly, that is,
to all orders in $\Delta t$. Hence, in this approach the generated
probability distribution differs from the stochastic distribution
by terms of order $\Delta t^{2}$, even in the infinite ensemble limit.

Hence, we use the parallel SDE to provide the \emph{initial} estimate
for the stochastic trajectory, prior to error optimization, of:
\begin{equation}
\mathbf{X}_{(0)}\left(t+\Delta t\right)=\mathbf{X}\left(t\right)+\mathbf{A}\left(t\right)\Delta t+\bm{\mathsf{B}}\left(t\right)\Delta\mathbf{W}\,.
\end{equation}

Obtaining the improved trajectory estimates is then essentially identical
to the static requirement given already. We can write the moment requirement
in the form~(\ref{eq:static-target-eqn}):
\begin{equation}
\bar{\mathbf{o}}\left(\mathbf{X}\left(t+\Delta t\right)\right)=\mathbf{c},\label{eq:combined-target-eqn}
\end{equation}
where $\mathbf{c}$ is a vector of ideal observables, and $\mathbf{X}\left(t+\Delta t\right)=\mathbf{X}_{(0)}\left(t+\Delta t\right)+\delta\mathbf{X}$.
Since we do not know these ``ideal'' averages over infinite number
of trajectories, we approximate $\mathbf{c}$ by assuming we know
all the moments exactly at the start of the time interval based on
knowing $\mathbf{X}(t)$. Changes in moments are estimated up to terms
of order $\Delta t$ using the approximations~(\ref{eq:ideal-observable-changes-appr-dt})
and~(\ref{eq:ideal-observable-changes-appr-N}) so that, in general:
\begin{eqnarray}
c_{m} & = & \left\langle o_{m}+\left[\frac{\partial o_{m}}{\partial x_{i}}a_{i}+\frac{1}{2}\frac{\partial^{2}o_{m}}{\partial x_{i}\partial x_{j}}\mathsf{d}_{ij}\right]\Delta t\right\rangle _{N_{S}}.
\end{eqnarray}

In the one-dimensional example case, from Eq.~(\ref{eq:ideal moment changes}),

\begin{eqnarray}
c_{m} & = & \left\langle x^{m}+m\left[x^{m-1}a+\frac{m-1}{2}x^{m-2}\mathsf{b}^{2}\right]\Delta t\right\rangle _{N_{S}}.\label{eq:combined:ideal-new-moments}
\end{eqnarray}

At the $i$-th step in time, we assume that $\mathbf{X}_{(i+1)}\left(t+\Delta t\right)=\mathbf{X}_{(i)}\left(t+\Delta t\right)+\delta\mathbf{X}_{(i+1)}$,
so that
\begin{equation}
\bar{\mathbf{o}}\left(\mathbf{X}_{(i+1)}\left(t+\Delta t\right)\right)=\bar{\mathbf{o}}\left(\mathbf{X}_{(i)}\left(t+\Delta t\right)\right)+\bm{\mathsf{J}}_{(i)}\delta\mathbf{X}_{(i+1)}+O(\delta\mathbf{X}_{(i+1)}^{2}).
\end{equation}

\subsection{Iterative solution to the combined equations}

The iterative equations to be solved are:

\begin{equation}
\bm{\mathsf{J}}_{(i)}\delta\mathbf{X}_{(i+1)}=\mathbf{c}-\bar{\mathbf{o}}_{(i)}.\label{eq:dynamic:combined-iteration}
\end{equation}
Here we define the matrices $\bm{\mathsf{J}}_{(i)}$ as in Eq.~(\ref{eq:XY-matrices}).
The moment equation is then solved iteratively using the techniques
outlined above, that is, either by pseudo-inverse iterations, or more
efficiently by the matrix inversion method. With pseudo-inverses,
the iteration equations are then simply:

\begin{equation}
\mathbf{X}_{(i+1)}=\mathbf{X}_{(i)}+\bm{\mathsf{J}}_{(i)}^{+}\left(\mathbf{c}-\bar{\mathbf{o}}_{(i)}\right).
\end{equation}
For cases where $\bm{\mathsf{u}}=\bm{\mathsf{J}}\bm{\mathsf{J}}^{\dagger}$
is invertible, we can define, just as in the static case:
\begin{equation}
\mathbf{X}_{(i+1)}=\mathbf{X}_{(i)}+\bm{\mathsf{J}}_{(i)}^{\dagger}\bm{\mathsf{u}}_{(i)}^{-1}\left(\mathbf{c}-\bar{\mathbf{o}}_{(i)}\right).
\end{equation}

Numerical results using this method are given in the next Section.

\subsection{Individual optimization\label{sub:dynamic:individual}}

In the individual optimization approach, both the finite ensemble
moments and the infinite ensemble moments are calculated to the same
order in $\Delta t$, which allows us to optimize each order in time
individually. This strategy permits a clear separation of ensemble
and time-step errors. For this purpose, it is convenient to define
an effective noise term, $\Delta\mathbf{V}=\bm{\mathsf{B}}\Delta\mathbf{W}$.
This is optimized to give the final change in $\mathbf{X}$.

From the moment equations of Eq.~(\ref{eq:moment-equations}) and~(\ref{eq:ideal-observable-changes-appr-dt}),
one has the required observable changes to order $\Delta t$:
\begin{equation}
\sum_{n}\frac{\partial\bar{o}_{m}}{\partial X_{n}}\Delta V_{n}\left(\mathbf{X}\right)+\frac{1}{2}\sum_{n,p}\frac{\partial^{2}\bar{o}_{m}}{\partial X_{n}\partial X_{p}}\left[\Delta V_{n}\left(\mathbf{X}\right)\Delta V_{p}\left(\mathbf{X}\right)-\Delta t\mathsf{D}_{np}\left(\mathbf{X}\right)\right]=e_{m}=\mathcal{O}\left(\Delta t^{2}\right)
\end{equation}

Here, since we wish to set every moment error to zero over a finite
set of moments, it is convenient to define an error vector, $\mathbf{e}=\left(e_{1},\ldots e_{M}\right)^{T}$.
For example, in the one-dimensional case, from Eq.~(\ref{eq:moment-equations})
the extended diffusion matrix is simply diagonal, and one has that:
\begin{equation}
\frac{m}{N_{S}}\sum_{n}\left[X_{n}^{m-1}\Delta V_{n}+\frac{m-1}{2}X_{n}^{m-2}\left(\Delta V_{n}^{2}-\mathsf{D}_{nn}\Delta t\right)\right]_{N_{S}}=e_{m}=\mathcal{O}\left(\Delta t^{2}\right),\label{eq:dynamic:moment-equations}
\end{equation}

The error terms can be viewed as two distinct ``mean'' and ``variance''
conditions analogous to those in Eq.~(\ref{eq:noise-averages}),
and are equal to zero in the limit of an infinite ensemble. The error
requirement is satisfied provided two individual conditions are met,
representing terms of order $\sqrt{\Delta t}$ and $\Delta t$ respectively:

\begin{eqnarray}
\mathbf{e}^{(1)} & = & \bm{\mathsf{J}}\Delta\mathbf{V}^{\mathrm{(opt)}}=\mathbf{0}\nonumber \\
\mathbf{e}^{(2)} & = & \frac{1}{2}\bm{\mathcal{H}}:\left[\Delta\mathbf{V}^{\mathrm{(opt)}}\left(\Delta\mathbf{V}^{\mathrm{(opt)}}\right)^{T}-\bm{\mathsf{D}}\Delta t\right]=\mathbf{0},\label{eq:dynamic:individual:target-eqn}
\end{eqnarray}
where ``:'' stands for double contraction, so $\left(\bm{\mathcal{H}}:\mathbf{V}\mathbf{V}^{T}\right)_{i}\equiv\sum_{jk}\mathcal{H}_{ijk}V_{j}V_{k}$.
Here, we define $\bm{\mathsf{J}}$ as previously, and $\bm{\mathcal{H}}$
as a Hessian tensor:
\begin{equation}
\mathcal{H}_{mnp}\left(\mathbf{X}\right)=\frac{\partial^{2}\bar{o}_{m}\left(\mathbf{X}\right)}{\partial X_{n}\partial X_{p}}.
\end{equation}
In the one-dimensional moment-based example, the results are simpler.
Due to the absence of cross-derivative terms, the second derivative
tensor becomes a matrix,
\begin{equation}
\mathsf{H}_{mn}\left(\mathbf{X}\right)=\frac{\partial^{2}\bar{o}_{m}\left(\mathbf{X}\right)}{\partial X_{n}^{2}}=m\frac{m-1}{N_{S}}X_{n}^{m-2}.
\end{equation}
The variance condition then becomes:

\begin{equation}
\mathbf{e}^{(2)}\equiv\frac{1}{2}\boldsymbol{\mathsf{H}}\left(\Delta\mathbf{V}^{\mathrm{(opt)}}\circ\Delta\mathbf{V}^{\mathrm{(opt)}}-\mathbf{D}\Delta t\right)=\mathbf{0},\label{eq:opt-conditions-matrix-individual-1}
\end{equation}
where $D_{n}\equiv\mathsf{D}_{nn}$.

\subsection{Iterative solution to the individual equations}

These are nonlinear equations, and to solve them we implement an iterative
approach. For a large ensemble, the optimizing equations are nearly
satisfied to zero-th order, up to errors of order $1/\sqrt{N_{S}}$.
Hence, it is efficient to iteratively solve the equations by linearization,
giving a Newton-Raphson approach similar to the combined method given
above.

Defining the difference $\delta\mathbf{X}_{(i+1)}\equiv\Delta\mathbf{V}_{(i+1)}-\Delta\mathbf{V}_{(i)}$,
the separate conditions are linearized by assuming that $\|\delta\mathbf{X}_{(i+1)}\|\ll\|\Delta\mathbf{V}\|$,
to give

\begin{eqnarray}
\bm{\mathsf{J}}\left(\Delta\mathbf{V}_{(i)}+\delta\mathbf{X}_{(i+1)}\right) & = & \boldsymbol{0},\nonumber \\
\frac{1}{2}\bm{\mathcal{H}}:\left(\Delta\mathbf{V}_{(i)}\Delta\mathbf{V}_{(i)}^{T}+\delta\mathbf{X}_{(i+1)}\Delta\mathbf{V}_{(i)}^{T}+\Delta\mathbf{V}_{(i)}\delta\mathbf{X}_{(i+1)}^{T}-\bm{\mathsf{D}}\Delta t\right) & = & \boldsymbol{0}.
\end{eqnarray}
Noting that $\boldsymbol{\mathcal{H}}$ is symmetric with respect
to its last two indices ($\mathcal{H}_{mnp}=\mathcal{H}_{mpn}$),
we can rewrite the terms with $\delta\mathbf{X}_{(i+1)}$ in the second
equation as

\begin{equation}
\frac{1}{2}\bm{\mathcal{H}}:\left(\delta\mathbf{X}_{(i+1)}\Delta\mathbf{V}_{(i)}^{T}+\Delta\mathbf{V}_{(i)}\delta\mathbf{X}_{(i+1)}^{T}\right)\equiv\left(\boldsymbol{\mathcal{H}}\Delta\mathbf{V}_{(i)}\right)\delta\mathbf{X}_{(i+1)},
\end{equation}
which allows us to concatenate the two systems into a single matrix
equation
\begin{equation}
\tilde{\bm{\mathsf{J}}}_{(i)}\delta\mathbf{X}_{(i+1)}=\tilde{\mathbf{R}}\left(\Delta\mathbf{V}_{(i)}\right).\label{eq:delta-matrix-expansion-1-1}
\end{equation}
Here the doubly extended matrix $\tilde{\bm{\mathsf{J}}}$ is:
\begin{equation}
\tilde{\bm{\mathsf{J}}}=\left(\begin{array}{c}
\bm{\mathsf{J}}\\
\bm{\mathcal{H}}\Delta\mathbf{V}
\end{array}\right).\label{eq:noise-transformation-matrix-1-1}
\end{equation}
and the remainder vector $\tilde{\mathbf{R}}\left(\Delta\mathbf{V}\right)$
is defined as:

\begin{equation}
\tilde{\mathbf{R}}\left(\Delta\mathbf{V}\right)=\left(\begin{array}{c}
-\bm{\mathsf{J}}\Delta\mathbf{V}\\
\frac{1}{2}\bm{\mathcal{H}}:\left(\mathbf{D}\Delta t-\Delta\mathbf{V}_{(i)}\Delta\mathbf{V}_{(i)}^{T}\right)
\end{array}\right).
\end{equation}

Starting with $i=0$, we now take $\Delta\mathbf{V}_{(i)}$ and solve
these equations iteratively to obtain $\Delta\mathbf{V}_{(i+1)}=\delta\mathbf{X}_{(i+1)}+\Delta\mathbf{V}_{(i)}$.
To obtain a linear solution for each iterative step, we can use a
pseudo-inverse as before, so that:

\begin{equation}
\Delta\mathbf{V}_{(i+1)}=\Delta\mathbf{V}_{(i)}+\tilde{\bm{\mathsf{J}}}_{(i)}^{+}\tilde{\mathbf{R}}\left(\Delta\mathbf{V}_{(i)}\right).
\end{equation}
As previously, this reduces to ordinary matrix inversion in most cases
of interest since, we can define a positive semi-definite matrix

\begin{eqnarray}
\tilde{\bm{\mathsf{u}}}\left(\Delta\mathbf{v}\right) & = & \tilde{\bm{\mathsf{J}}}\left(\Delta\mathbf{V}\right)\tilde{\bm{\mathsf{J}}}^{\dagger}\left(\Delta\mathbf{V}\right).
\end{eqnarray}
Provided this is invertible, this gives an iterative procedure for
the individual POS algorithm, which is:

\begin{equation}
\Delta\mathbf{V}_{(i+1)}=\Delta\mathbf{V}_{(i)}+\tilde{\bm{\mathsf{J}}}_{(i)}^{\dagger}\tilde{\bm{\mathsf{u}}}_{(i)}^{-1}\mathbf{R}_{(i)}.
\end{equation}
This method can be implemented to scale as $\mathcal{O}\left(MN_{S}\right)$
by pre-calculating repeating matrix elements.

\section{Synthetic one-step benchmarks\label{sec:Synthetic-one-step-benchmarks}}

Before we apply our POS method to solving an SDE with error-propagation
over many successive minimization steps, it is informative to see
how effective it is at finding a solution for the corresponding single
step equation.

In order to evaluate the statistical performance of these two methods
at single-step error optimization, we implement the combined and individual
optimization methods with $10000$ separate initial ensembles, each
of $1000$ parallel samples of $\mathbf{X}$ and initial stochastic
noises $\mathbf{W}$. In all the tests, the initial elements of $\mathbf{X}$
are normally distributed with the mean $1$ and standard deviation
$0.1$, $a_{i}=0.5$ for all $i$ (constant drift), $\mathsf{b}_{i}=0.5$
for all $i$ (additive noise).

Results for more complicated synthetic tests show similar behavior,
so we do not include them here.

\subsection{Combined optimization}

In this section we are interested in how well we can solve the target
equation~(\ref{eq:combined-target-eqn}). The error measure in this
case is essentially the same as in Section~\ref{sub:benchmarks:static}:

\begin{equation}
\tilde{R}_{m}\left(\mathbf{X}\right)=\left|\frac{1}{N_{S}}\mathbf{X}^{m}\cdot\mathbf{1}-c_{m}\right|,
\end{equation}
where the target moments $c_{m}$ are given by Eq.~(\ref{eq:combined:ideal-new-moments}).

\begin{figure}
\includegraphics{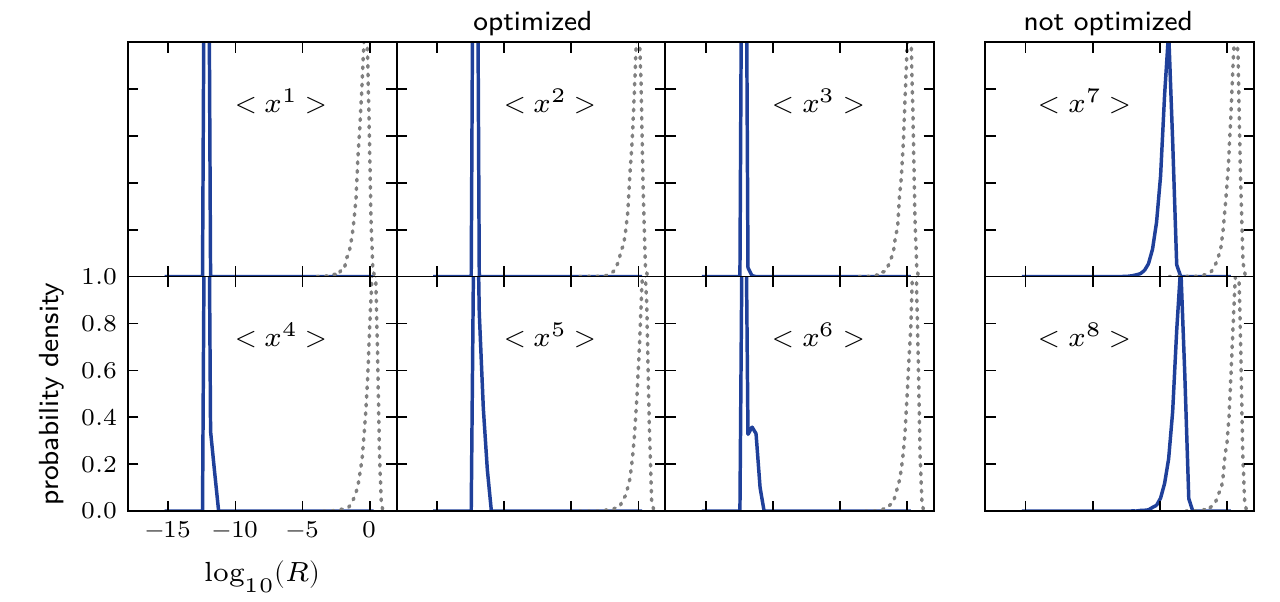}

\protect\caption{\label{fig:synthetic:combined:errors}Distribution of normalized errors
in the synthetic benchmark for the combined method with $\Delta t=10^{-4}$,
$N_{S}=1000$. Blue solid lines show the distribution of $R_{m}$
for the optimized vector, grey dotted lines show the distribution
of $R_{m}$ for the initial approximation $\mathbf{X}_{(0)}\left(t+\Delta t\right)$. }
\end{figure}

Similarly to the previous section, we will scale the errors as

\begin{equation}
R_{m}=\frac{\tilde{R}_{m}}{\sqrt{\Delta t/N_{S}}}
\end{equation}
in order to bring the errors for the unoptimized, randomly sampled
vectors close to $1$.

The results in Fig.~\ref{fig:synthetic:combined:errors} show that
a significant improvement is achieved for all $6$ moments being optimized,
with normalized sampling errors reduced to $10^{-12}$, or around
twelve orders of magnitude in these examples. Even the non-optimized
moments of orders $m=7,8$ have errors reduced to $10^{-5}$, or around
five orders.

\subsection{Number of iterations and convergence distance}

\begin{figure}
\includegraphics{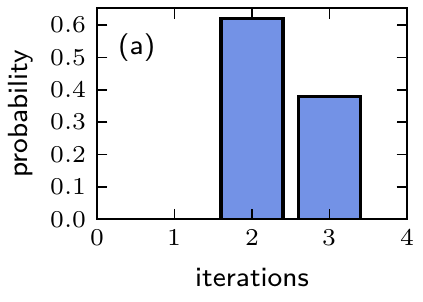}\includegraphics{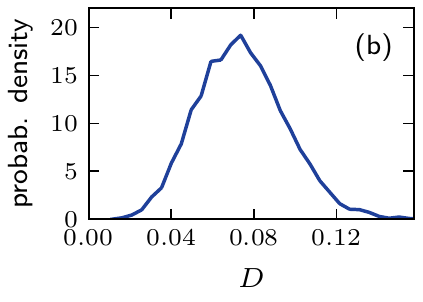}\includegraphics{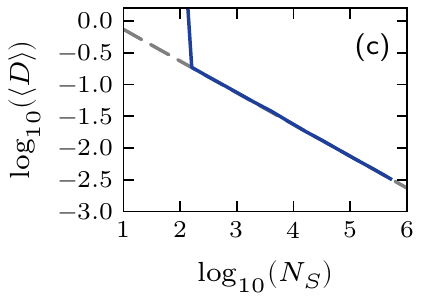}

\protect\caption{\label{fig:synthetic:combined:distance-and-iterations}(a) number
of iterations taken in the synthetic test for the combined method,
$\Delta t=10^{-4}$, $N_{S}=1000$, (b) Relative distance $D$ between
the optimized vector $\mathbf{X}^{\mathrm{(opt)}}\left(t+\Delta t\right)$
and the initial approximation $\mathbf{X}_{(0)}\left(t+\Delta t\right)$,
and (c) average relative distance $\langle D\rangle$ depending on
the number of trajectories $N_{S}$ in the synthetic optimization
tests for the individual expansion method, $\Delta t=10^{-4}$. The
grey dashed line denotes the slope of a $1/\sqrt{N_{S}}$ dependence.}
\end{figure}

The method requires only a few iterations to converge, as Fig.~\ref{fig:synthetic:combined:distance-and-iterations}(a)
demonstrates.

Next, we calculate the distance of the variance-reduced solution from
the original estimate. It is important that this quantity is as small
as possible, consistent with the targeted variance reduction, to eliminate
systematic errors in higher-order, non-minimized moments. We compute:

\begin{equation}
D=\frac{\|\Delta\mathbf{X}^{\mathrm{(opt)}}-\Delta\mathbf{X}_{(0)}\|}{\|\Delta\mathbf{X}_{(0)}\|}\equiv\frac{\|\mathbf{X}^{\mathrm{(opt)}}\left(t+\Delta t\right)-\mathbf{X}_{(0)}\left(t+\Delta t\right)\|}{\|\mathbf{X}_{(0)}\left(t+\Delta t\right)-\mathbf{X}\left(t\right)\|}
\end{equation}
to estimate how far the solution is from the initial approximation.
Fig.~\ref{fig:synthetic:combined:distance-and-iterations}(b) shows
that the deviation of the optimized noise from the original is $0.07$
on average, which is of order $1/\sqrt{N_{s}}$ as expected.

\begin{figure}
\includegraphics{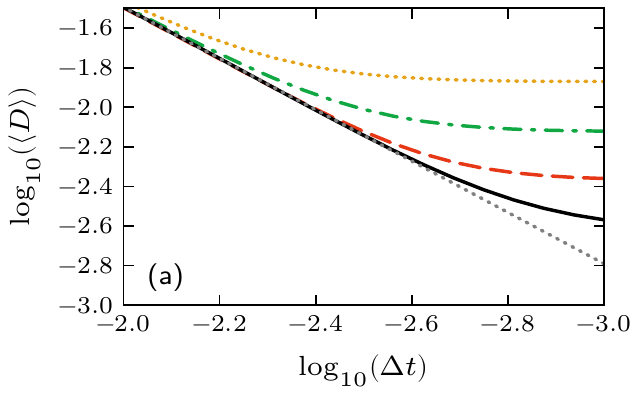}\includegraphics{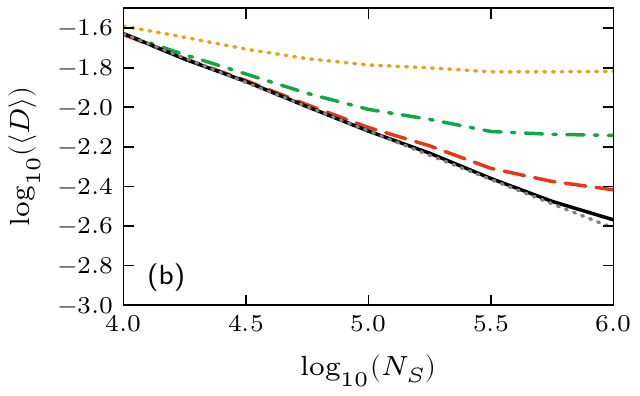}

\protect\caption{\label{fig:synthetic:combined:distance-plateau}(a) Average relative
distance $\langle D\rangle$ depending on $\Delta t$ in the synthetic
optimization tests for the combined method, with $N_{S}=10^{4.5}$
(dotted yellow line), $N_{S}=10^{5}$ (dash-dotted green line), $N_{S}=10^{5.5}$
(dashed red line) and $N_{S}=10^{6}$ (solid black line). The grey
dashed line denotes the least squares fit for the linear part of the
$N_{S}=10^{6}$ line and corresponds to the dependence $\langle D\rangle\propto\Delta t^{1.29}$.
(b) Average relative distance $\langle D\rangle$ depending on $N_{S}$
in the synthetic optimization tests for the combined method, with
$\Delta t=10^{-2.25}$ (dotted yellow line), $\Delta t=10^{-2.5}$
(dash-dotted green line), $\Delta t=10^{-2.75}$ (dashed red line)
and $\Delta t=10^{-3}$ (solid black line). The grey dashed line denotes
the least squares fit for the linear part of the $\Delta t=10^{-3}$
line and corresponds to the dependence $\langle D\rangle\propto N_{S}^{-0.49}$.}
\end{figure}

The difference in the orders of $\Delta t$ in the left and the right
parts of the target equation~(\ref{eq:combined-target-eqn}) leads
to the distance $D$ being bounded by a finite value, depending on
either $\Delta t$ or $1/N_{S}$. This is unlike the behavior demonstrated
by the static optimization in the previous section. Our tests show
that for very small $\Delta t$, the asymptotic dependence of $\langle D\rangle$
on $N_{S}$ is $\langle D\rangle\propto1/\sqrt{N_{S}}$, as shown
in Fig.~\ref{fig:synthetic:combined:distance-plateau}(a). Similarly,
for very large $N_{S}$ the dependence is $\langle D\rangle\propto\Delta t^{3/2}$,
as seen in Fig.~\ref{fig:synthetic:combined:distance-plateau}(b).

This asymptotic error is essentially due to the fact that the moment
equations are themselves an expansion in $\Delta t$, leading to a
residual truncation error.

\begin{figure}
\includegraphics{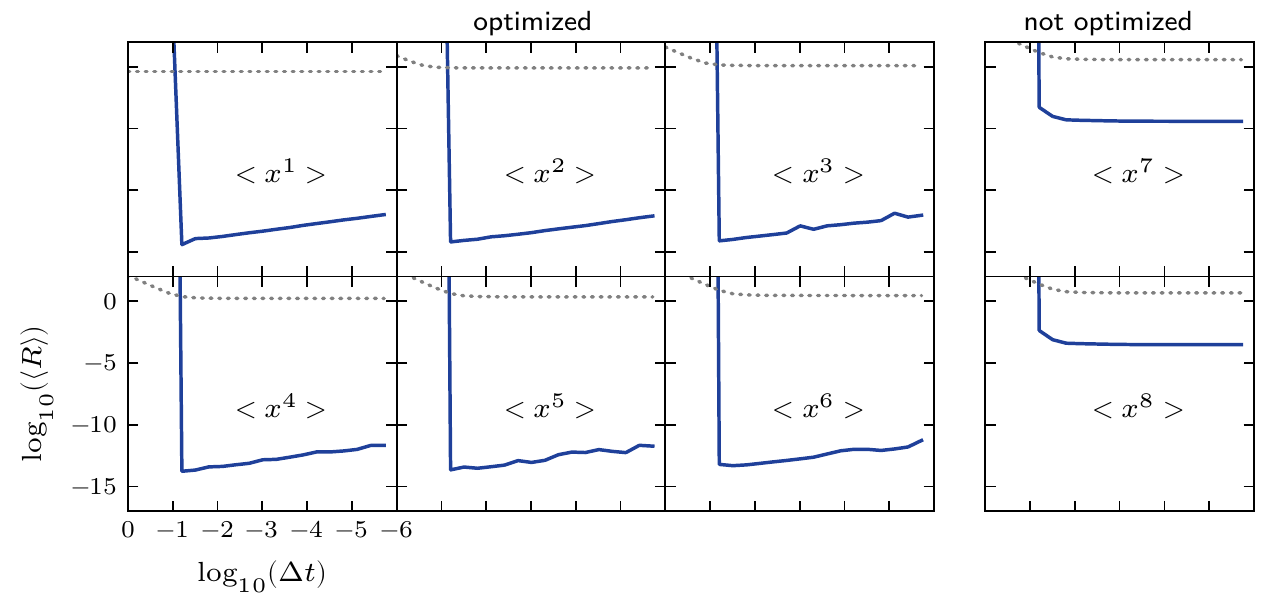}

\protect\caption{\label{fig:synthetic:combined:error-vs-dt}Dependence of errors on
the number of trajectories $N_{S}$ in the synthetic test for the
combined method, $\Delta t=10^{-4}$. Blue solid lines show the average
of $R_{m}$ for the optimized vector, grey dotted lines show the average
of $R_{m}$ for the initial approximation $\mathbf{X}_{(0)}\left(t+\Delta t\right)$.}
\end{figure}

\begin{figure}
\includegraphics{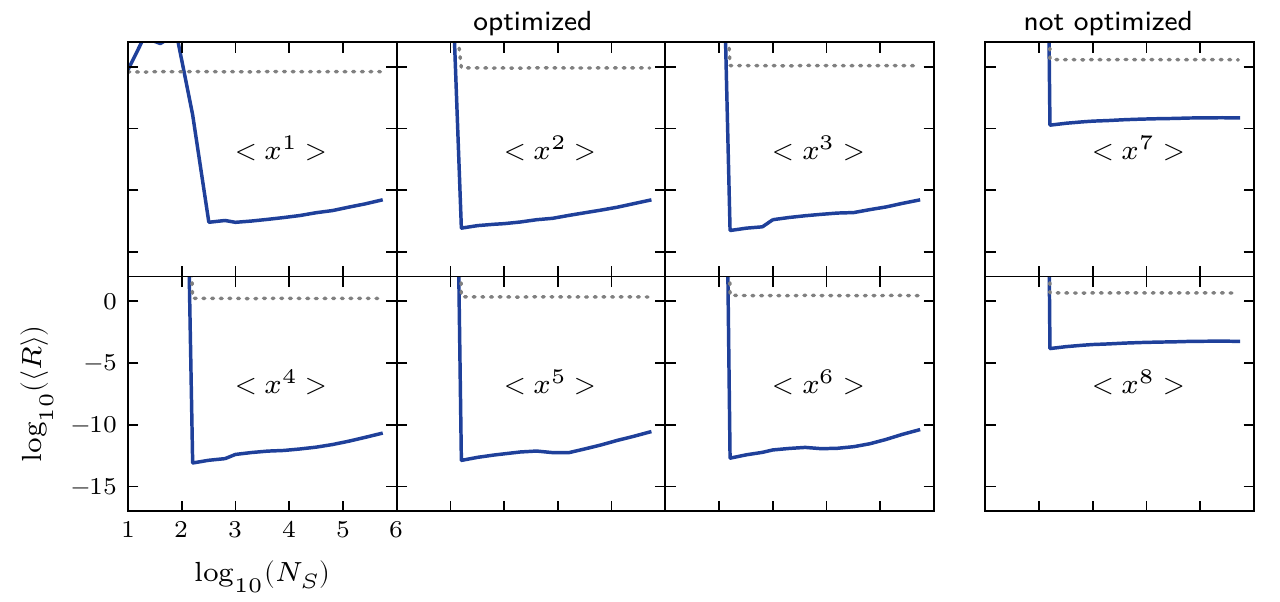}

\protect\caption{\label{fig:synthetic:combined:error-vs-N}Dependence of average scaled
errors on $\Delta t$ in the synthetic test for the combined method,
$N_{S}=1000$. Blue solid lines show the average of $R_{m}$ for the
optimized vector, grey dotted lines show the average of $R_{m}$ for
the initial approximation $\mathbf{X}_{(0)}\left(t+\Delta t\right)$.}
\end{figure}

The tests of the dependence on the noise scale ($\Delta t$) show
that the method is stable for a wide range of values, but breaks down
at larger values of $\mathrm{\Delta}t$ as seen in Fig.~\ref{fig:synthetic:combined:error-vs-dt}.
The dependence on the number of trajectories in Fig.~\ref{fig:synthetic:combined:error-vs-N}
displays a similar cutoff as seen in the tests of the static POS\@,
in the results of the previous subsection.

\subsection{Individual optimization}

In this subsection we will use the individual method from Section~\ref{sub:dynamic:individual},
for a similar type of synthetic benchmark. Since the target equations
we try to satisfy are different from the previous section, the definition
of the error changes to simply

\begin{equation}
\tilde{R}_{m}\left(\mathbf{X}\right)=|e_{m}^{(1)}|+|e_{m}^{(2)}|,
\end{equation}
where the vectors $\mathbf{e}^{(1)}$ and $\mathbf{e}^{(2)}$ are
defined by Eq.~(\ref{eq:dynamic:individual:target-eqn}).

Similarly to the approach of the previous subsection, we normalize
the error as

\begin{equation}
R_{m}=\frac{\tilde{R}_{m}}{\sqrt{\Delta t/N_{S}}}.
\end{equation}

\begin{figure}
\includegraphics{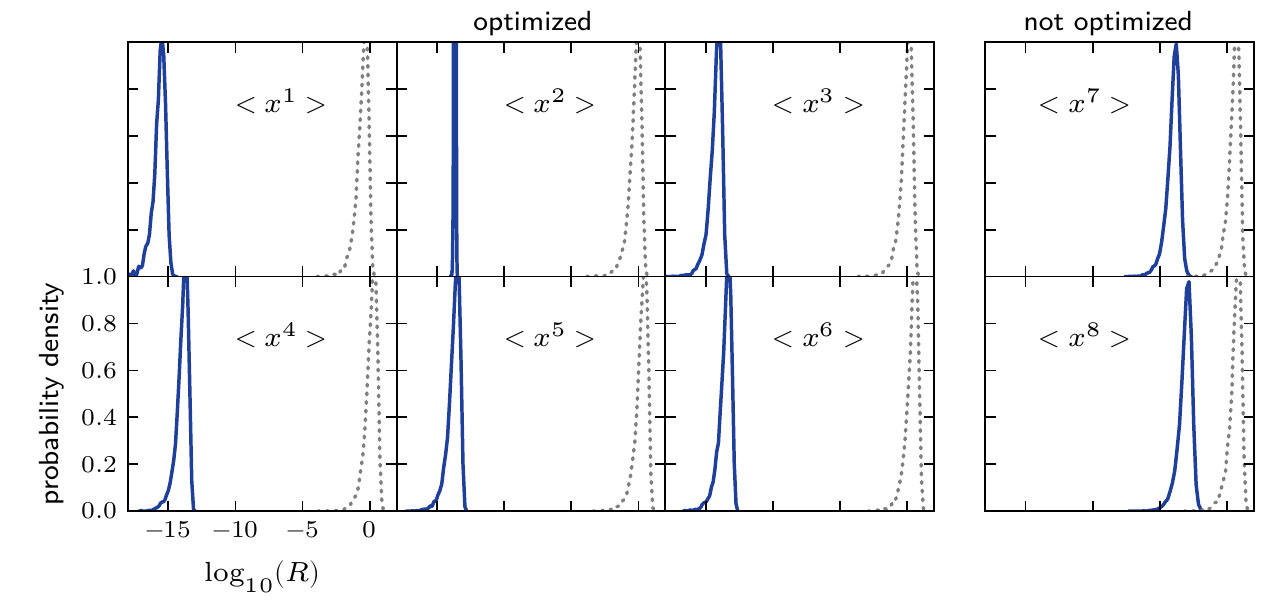}

\protect\caption{\label{fig:synthetic:individual:errors}Distribution of normalized
errors in the synthetic benchmark for the individual method with $\Delta t=0.1$,
$N_{S}=1000$. Blue solid lines show the distribution of $R_{m}$
for the optimized vector, grey dotted lines show the distribution
of $R_{m}$ for the initial effective noise term $\Delta\mathbf{V}$.}
\end{figure}

The results in Fig.~\ref{fig:synthetic:individual:errors}(a) show
that a significant improvement is achieved for the first $6$ moments
being optimized, with results only limited by the numerical precision.
To check that numerical precision was the limiting factor, tests were
run with the same method in quadruple precision. This demonstrated
that the improvement was indeed limited only by the numerical precision,
with errors now as low as $10^{-30}$. The stopping condition must
be changed to $\eta=10^{-15}$ to reach this greater numerical precision
limit.

\subsection{Number of iterations and convergence distance}

\begin{figure}
\includegraphics{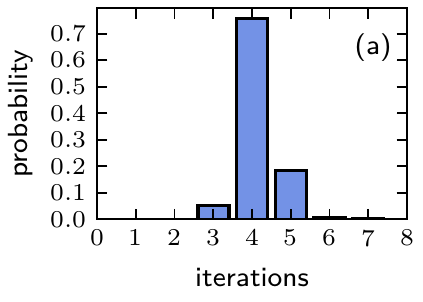}\includegraphics{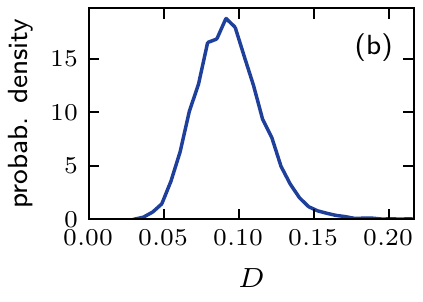}\includegraphics{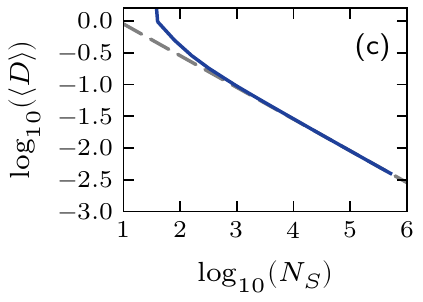}

\protect\caption{\label{fig:synthetic:individual:distance-and-iterations}(a) The number
of iterations taken in the synthetic test for the individual method,
$\Delta t=0.1$, $N_{S}=1000$. (b) Relative distance $D$ between
the optimized vector $\Delta\mathbf{V}^{\mathrm{(opt)}}$ and the
initial approximation $\Delta\mathbf{V}$. (c) Average relative distance
$\langle D\rangle$ depending on the number of trajectories $N_{S}$
in the synthetic optimization tests for the individual method, $\Delta t=0.1$.
The grey dashed line denotes the slope of a $1/\sqrt{N_{S}}$ dependence.}
\end{figure}

We also calculate the distance

\begin{equation}
D\equiv\frac{\|\Delta\mathbf{V}^{(\mathrm{opt})}-\Delta\mathbf{V}\|}{\|\Delta\mathbf{V}\|}
\end{equation}
and the number of iterations in our tests. The method also requires
only a few iterations to converge, as Fig.~\ref{fig:synthetic:individual:distance-and-iterations}(a)
demonstrates. Fig.~\ref{fig:synthetic:individual:distance-and-iterations}(b)
shows that the deviation of the optimized noise from the original
one is small. Unlike the combined method, there is no mismatch between
the orders of $\Delta t$ in Eq.~(\ref{eq:dynamic:individual:target-eqn}),
and the distance therefore converges as $\langle D\rangle\propto1/\sqrt{N_{S}}$,
as shown in Fig.~\ref{fig:synthetic:individual:distance-and-iterations}(c).

\begin{figure}
\includegraphics{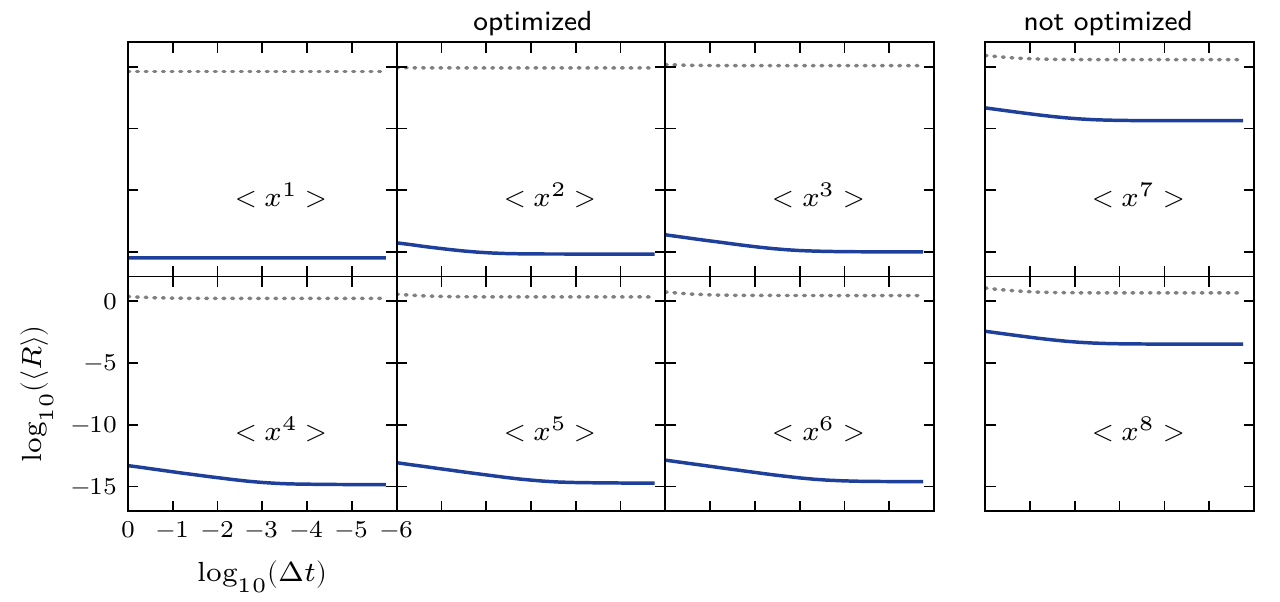}

\protect\caption{\label{fig:synthetic:individual:error-vs-dt}Dependence of average
scaled errors on $\Delta t$ in the synthetic test for the individual
method, $N_{S}=1000$. Blue solid lines show the average of $R_{m}$
for the optimized vector, grey dotted lines show the average of $R_{m}$
for the initial approximation $\Delta\mathbf{V}_{(0)}$.}
\end{figure}

\begin{figure}
\includegraphics{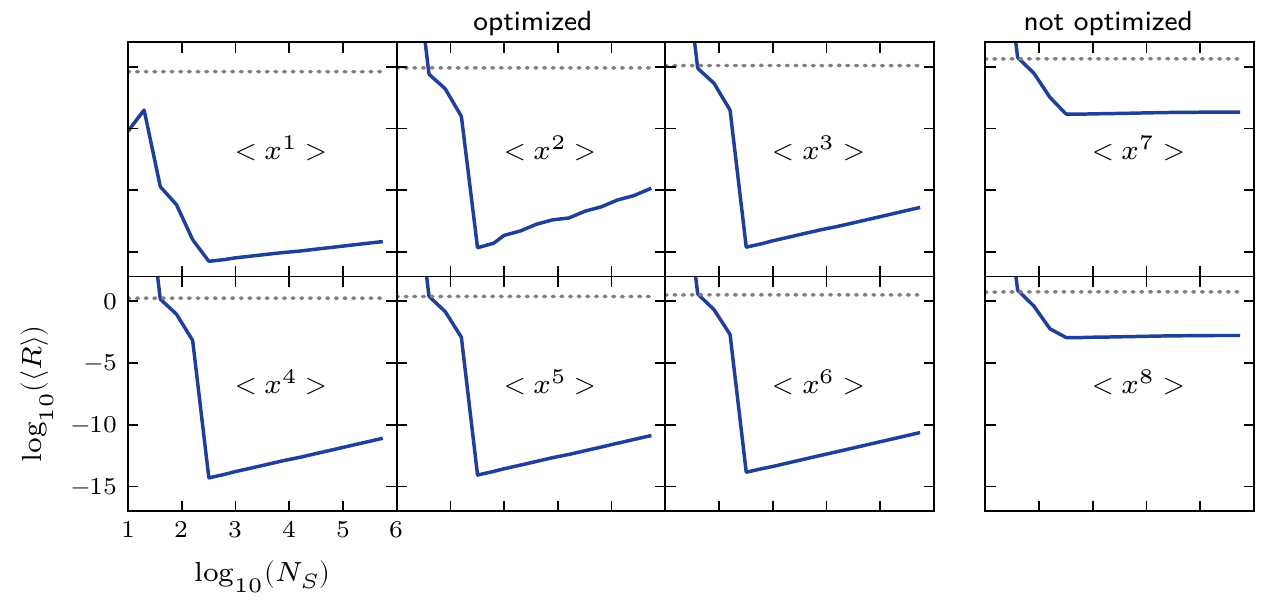}

\protect\caption{\label{fig:synthetic:individual:error-vs-N}Dependence of errors on
the number of trajectories $N_{S}$ in the synthetic test for the
expansion method, $\Delta t=0.1$. Blue solid lines show the average
of $R_{m}$ for the optimized vector, grey dotted lines show the average
of $R_{m}$ for the initial approximation $\Delta\mathbf{V}_{(0)}$.}
\end{figure}

The tests of the dependence on the noise scale ($\Delta t$) show
that this method is exceptionally stable over a wide range of step-size
values, as seen in Fig.~\ref{fig:synthetic:individual:error-vs-dt}.
The dependence on the number of trajectories in Fig.~\ref{fig:synthetic:individual:error-vs-N}
displays a similar cutoff to those seen in the tests of the static
POS\@.

\section{Dynamical numerical results\label{sec:benchmarks}}

Next, putting all the preliminary results together, tests were performed
of complete stochastic differential equation solutions. Although all
results demonstrate strong variance reduction, the amount and nature
of variance reduction depends on the precise equation, as explained
below.

\subsection{Linear drift, additive noise}

It is instructive to first consider an ``ideal'' linear SDE for
the application of POS, with the target observables being the moments
of $x$ as previously. An ideal linear case is one with a linear $a$
and $b$. This is preferred as a first test because the errors from
non-optimized moments do not get mixed with the optimized ones, as
explained in Section~\ref{sub:dynamic:error-propagation}. In these
benchmarks, both the combined and individual methods were used, giving
very similar performance overall.

As an example, we used a one-dimensional Ornstein-Uhlenbeck process

\begin{equation}
\mathrm{d}x=\left(f-gx\right)\mathrm{d}t+b\mathrm{d}w,
\end{equation}
where $\langle\mathrm{d}w\rangle=0$ and $\langle\mathrm{d}w^{2}\rangle=dt$,
and $f$, $g$ and $b$ are constants. It has well-known exact solutions
for every moment. Namely, for Gaussian starting conditions the mean
is given by:

\begin{equation}
\langle x\rangle=\frac{f}{g}\left(1-e^{-gt}\right)+e^{-gt}\langle x\rangle_{t=0},
\end{equation}
and for the variance,

\begin{equation}
\langle\left(x-\langle x\rangle\right)^{2}\rangle=e^{-2gt}\langle\left(x-\langle x\rangle\right)^{2}\rangle_{t=0}+\frac{b^{2}}{2g}\left(1-e^{-2gt}\right),
\end{equation}

Higher order central moments are then calculated according to

\begin{equation}
\langle\left(x-\langle x\rangle\right)^{m}\rangle=0,\quad m=3,5,\dots
\end{equation}
and:

\begin{equation}
\langle\left(x-\langle x\rangle\right)^{m}\rangle=m!!\langle\left(x-\langle x\rangle\right)^{2}\rangle^{m/2},\quad m=4,6,\dots
\end{equation}
Consequently, the raw moments are calculated from the central ones
using the recursive formula

\begin{equation}
\langle x^{m}\rangle=\langle\left(x-\langle x\rangle\right)^{m}\rangle-\sum_{p=0}^{m-1}\left(\begin{array}{c}
m\\
p
\end{array}\right)\langle x^{p}\rangle\left(-\langle x\rangle\right)^{m-p}.
\end{equation}

We also use cumulants as benchmarks, since these are a commonly used
alternative statistical measure to the ordinary moment. The cumulants
are calculated from the moments recursively, as

\begin{equation}
\kappa_{m}=\langle x^{m}\rangle-\sum_{p=1}^{m-1}\left(\begin{array}{c}
m-1\\
p-1
\end{array}\right)\kappa_{p}\langle x^{m-p}\rangle.
\end{equation}

In our tests we took $f=1$, $g=0.2$, $b=0.5$. At initial times,
$x(0)$ is normally distributed with mean $0.5$ and standard deviation
$0.1$. We pre-optimized the starting distributions using the static
method from Section~\ref{sub:static-pos:expansion} in order to satisfy
the requirements of POS methods.

\subsection{Cumulant and moment variance reduction}

For our tests, we compare the deviations~\textemdash{} due to sampling
errors~\textemdash{} from the exact solution in the first $8$ moments
and first $8$ cumulants of $x$ at $t=1$ for the POS integrator.
We use both combined and individual methods, with a reference SDE
solver having the same number of samples, for comparison. The latter
uses the same integration method (explicit Euler) and the same number
of trajectories as POS\@. We set up POS to optimize the first $6$
moments, so that we could look at the behavior of both optimized and
non-optimized moments.

\begin{figure}
\includegraphics{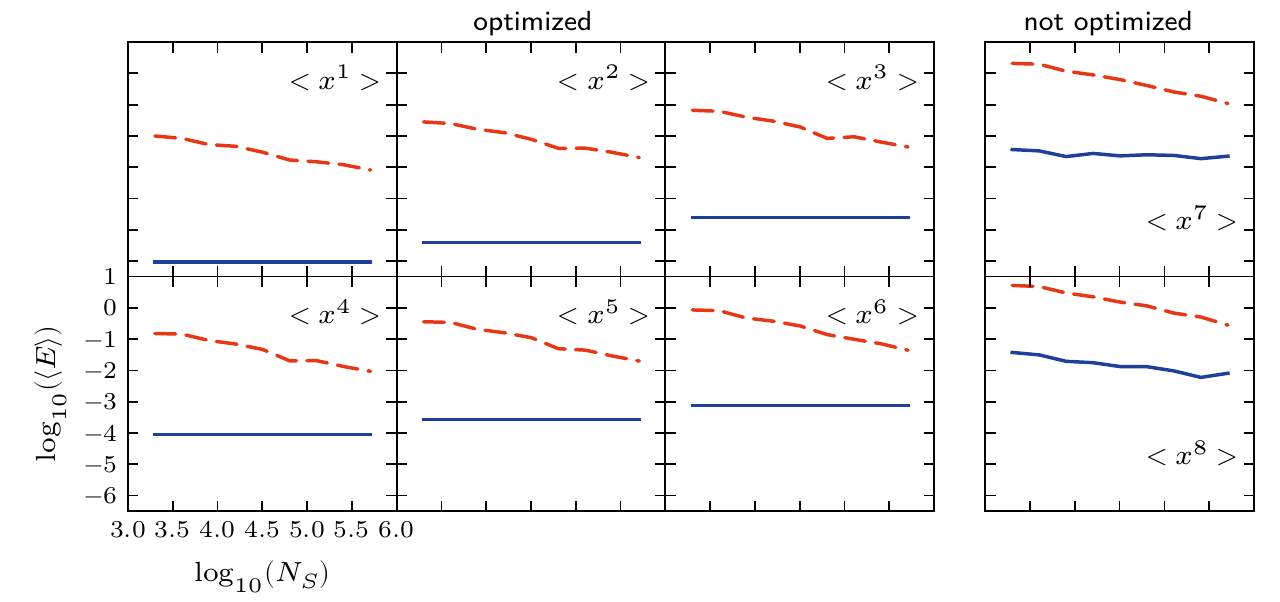}

\protect\caption{\label{fig:sde:linear:combined:diffs-moments}Linear SDE case with
$M=6$ optimized moments, showing the difference of the first $8$
moments from the exact solution at $t=1$ for the POS integrator using
the combined method (solid blue lines) and a reference explicit Euler
integrator (dashed red lines), with identical numbers of trajectories
and samples, averaged over $8$ tests. Both integrators use $N_{T}=8\times10^{4}$
time steps.}
\end{figure}

\begin{figure}
\includegraphics{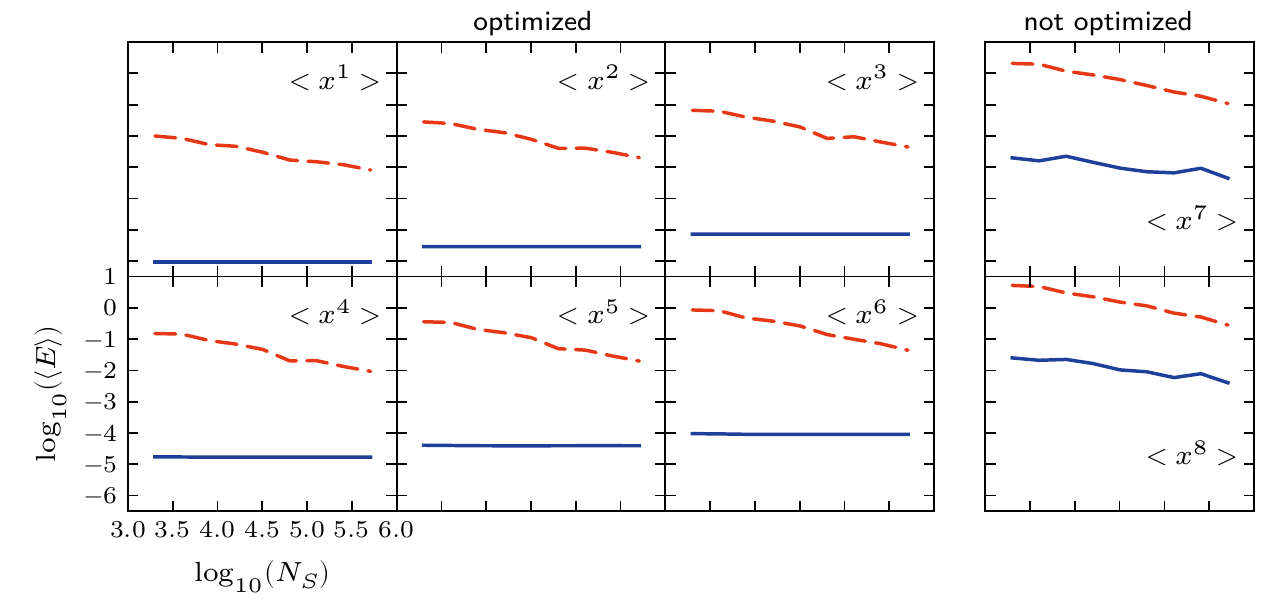}

\protect\caption{\label{fig:sde:linear:individual:diffs-moments}Linear SDE case with
$M=6$ optimized moments, showing the difference of the first $8$
moments from the exact solution at $t=1$ for the POS integrator using
the individual method (solid blue lines) and a reference explicit
Euler integrator (dashed red lines), with identical numbers of trajectories
and samples, averaged over $8$ tests. Both integrators use $N_{T}=8\times10^{4}$
time steps.}
\end{figure}

In Figs.~\ref{fig:sde:linear:combined:diffs-moments} and~\ref{fig:sde:linear:individual:diffs-moments}
we plotted the error in the first moments

\begin{equation}
E\left[\langle x^{m}\rangle\right]=\left|\langle x^{m}\rangle-\langle x^{m}\rangle^{(\mathrm{exact})}\right|
\end{equation}
for the reference and the two POS integrators, averaged over $8$
independent tests. In case of the POS methods the error is dominated
by the time step error $\epsilon_{T}$, while the error for the reference
integrator decreases proportional to $1/\sqrt{N_{S}}$. One can see
that the POS error is significantly decreased even for the moments
that were not directly optimized. The error improvement is largest
for the individual method, but is still several orders of magnitude
even for the simpler combined method. For the remaining graphs we
will mostly focus on the individual method which give better results,
to minimize the number of graphs, as both have similar general behavior.

\begin{figure}
\includegraphics{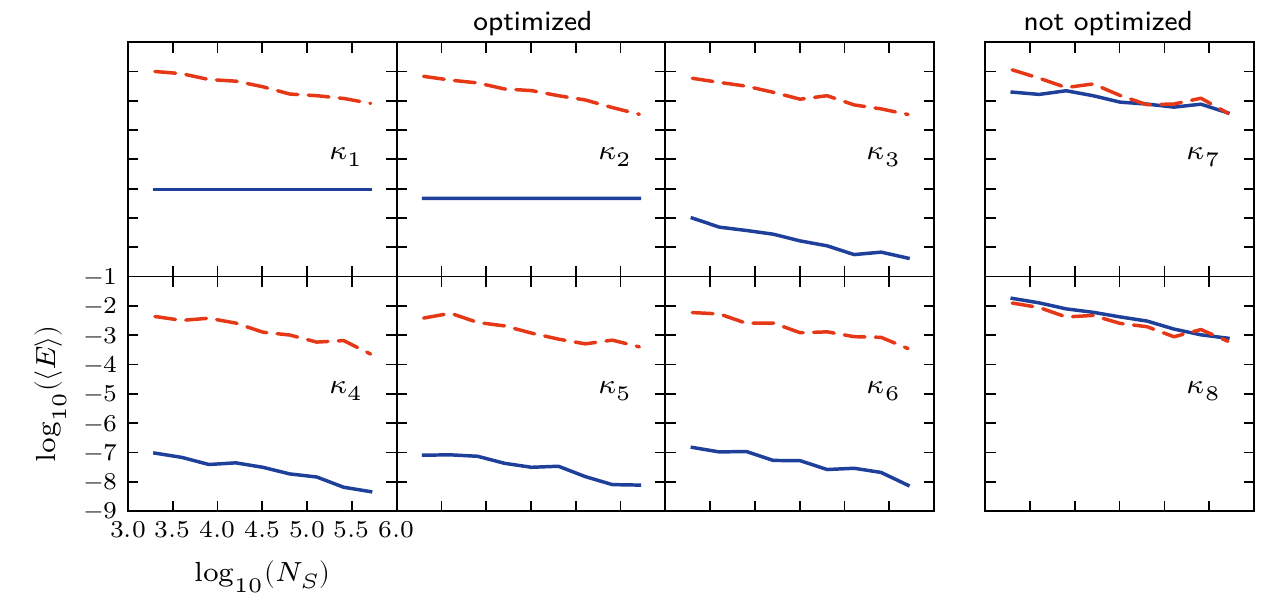}

\protect\caption{\label{fig:sde:linear:individual:diffs-cumulants}Linear SDE case
with $M=6$ optimized moments, showing the difference of the first
$8$ cumulants from the exact solution at $t=1$ for the POS integrator
using the individual method (solid blue lines) and a reference explicit
Euler integrator (dashed red lines), with identical numbers of trajectories
and samples, averaged over $8$ tests. Both integrators use $N_{T}=8\times10^{4}$
time steps.}
\end{figure}

The general behavior in the linear case can be understood better if
we plot the errors in the cumulants $E\left[\kappa_{p}\right]$ instead,
as in Fig.~\ref{fig:sde:linear:individual:diffs-cumulants}. It shows
that the error in the first $6$ cumulants is greatly reduced, while
the non-optimized cumulants stay roughly the same. In Fig.~\ref{fig:sde:linear:individual:diffs-moments},
$7$-th and $8$-th moments, while not being optimized directly, still
depend on the lower cumulants, and thus benefit from their reduced
error. This behavior is caused by the fact that in a linear SDE, cumulants
of different orders are not coupled with each other; however moments
\emph{are} coupled, leading to the observed error performance.

\subsection{Nonlinear drift, additive noise\label{sub:Nonlinear-drift,-additive-noise}}

As an example of a real-world, nonlinear SDE we can apply the POS
method to, we take an equation with a nonlinear drift term and additive
noise:

\begin{equation}
\mathrm{d}x=\left(x-x^{3}\right)\mathrm{d}t+\mathrm{d}w,\label{eq:nonlin_SDE}
\end{equation}
where $\langle\mathrm{d}w\rangle=0$ and $\langle\mathrm{d}w^{2}\rangle=1$.
This is similar to a known case where moment hierarchy methods fail
to give correct results~\cite{Morillo2014-checking}.

For large times $t$, Eq.~(\ref{eq:nonlin_SDE}) converges to a steady-state
solution. For any given observable $\left\langle f\left(x\right)\right\rangle $,
we can find its steady-state value using the formula

\begin{equation}
\left\langle f\left(x\right)\right\rangle _{ss}=\frac{\int_{-\infty}^{\infty}f\left(x\right)e^{x^{2}-\frac{1}{2}x^{4}}\mathrm{d}x}{\int_{-\infty}^{\infty}e^{x^{2}-\frac{1}{2}x^{4}}\mathrm{d}x}.
\end{equation}

We initialize our ensemble with normally-distributed numbers with
mean $0$ and standard deviation $\frac{1}{\sqrt{2}}$ and integrate
Eq.~(\ref{eq:nonlin_SDE}) using a combined POS-method as well as
an unoptimized Ito-Euler method until $t=25$ and then compare with
the steady-state values. The POS-optimized integration includes an
initial static optimization. Integrating until $t=25$ ensures that
Eq.~(\ref{eq:nonlin_SDE}) predicts observables sufficiently close
to the steady-state values. We assume that after some integration
time, the time evolution of the observables will be dominated by a
transient which is exponentially converging to the steady-state solution.
A least-squares fit of the time-evolution of the $8$'th moment is
in excellent agreement with the time-evolution resulting from our
simulations. From the fitting parameters, we can read off a decay
of more that $50$ mean lifetimes at $t=25$. We conclude that the
deviation from the steady-state value due to finite integration time
is negligible compared to sampling errors. As in the previous section,
we plot the first $8$ moments and first $8$ cumulants of $x$ at
$t=25$, but we optimize only the first $6$ moments when the POS
is used. In order to show the behaviour of non-polynomial observables,
we also plot the exponential and absolute value of $x$.

\begin{figure}
\includegraphics{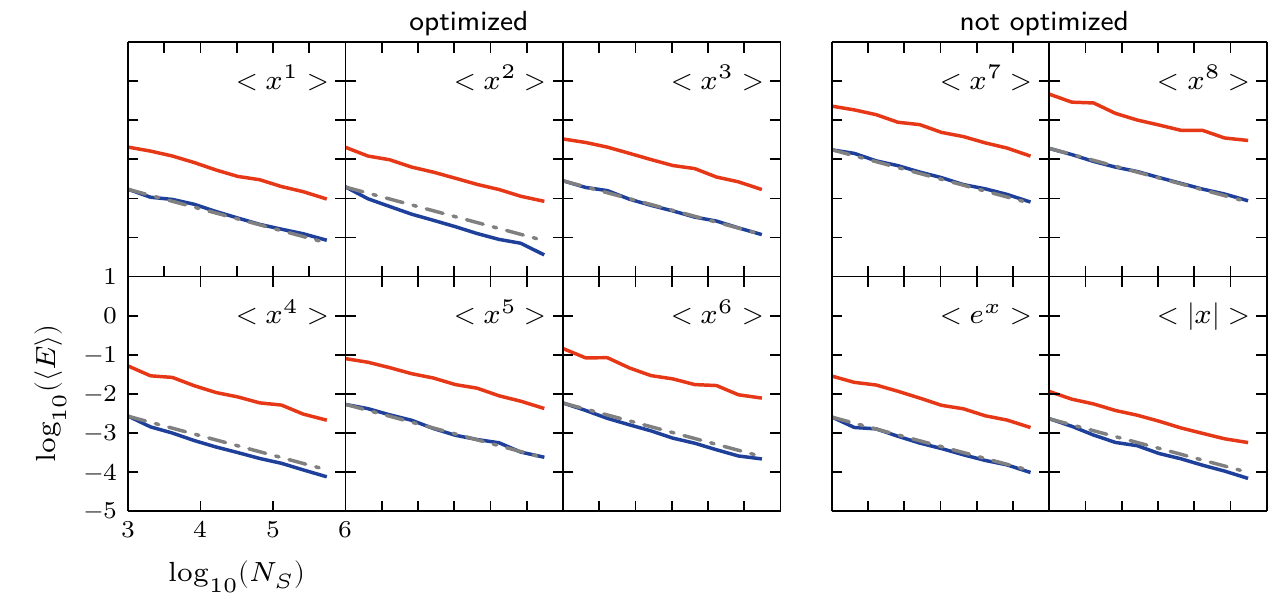}

\protect\caption{\label{fig:sde:nonlinear:combined:diffs-moments}Nonlinear SDE case
with $M=6$ optimized moments, showing the difference of the first
$8$ moments, as well as the exponential and absolute value function
from the steady-state solution at $t=25$ for the POS integrator using
the combined method (solid blue lines). We also include the errors
in an explicit Euler integrator (solid red lines), with identical
numbers of trajectories and samples, for comparison. Both integrators
use $N_{T}=12.5\times10^{3}$ time steps. The dash-dotted grey line
shows the slope for $\propto N_{S}^{-1/2}$. }
\end{figure}

In Fig.~\ref{fig:sde:nonlinear:combined:diffs-moments} we are plotting
differences from the accurate solution for the first $8$ moments:

\begin{equation}
E\left[\langle x^{m}\rangle\right]=\left|\langle x^{m}\rangle-\langle x^{m}\rangle_{ss}\right|\label{eq:E_xm}
\end{equation}
as well as for the exponential and absolute value function defined
analogously to Eq.~(\ref{eq:E_xm}) averaged over $120$ independent
test runs for the POS integration and $600$ runs for the reference
integration.

\textbf{The POS integrator consistently outperforms the reference
integrator, with errors reduced by a factor of $10-100$. }

A sampling error reduction of this order of magnitude requires $10^{2}-10^{4}$
times more samples using standard integrators.

\begin{figure}
\includegraphics{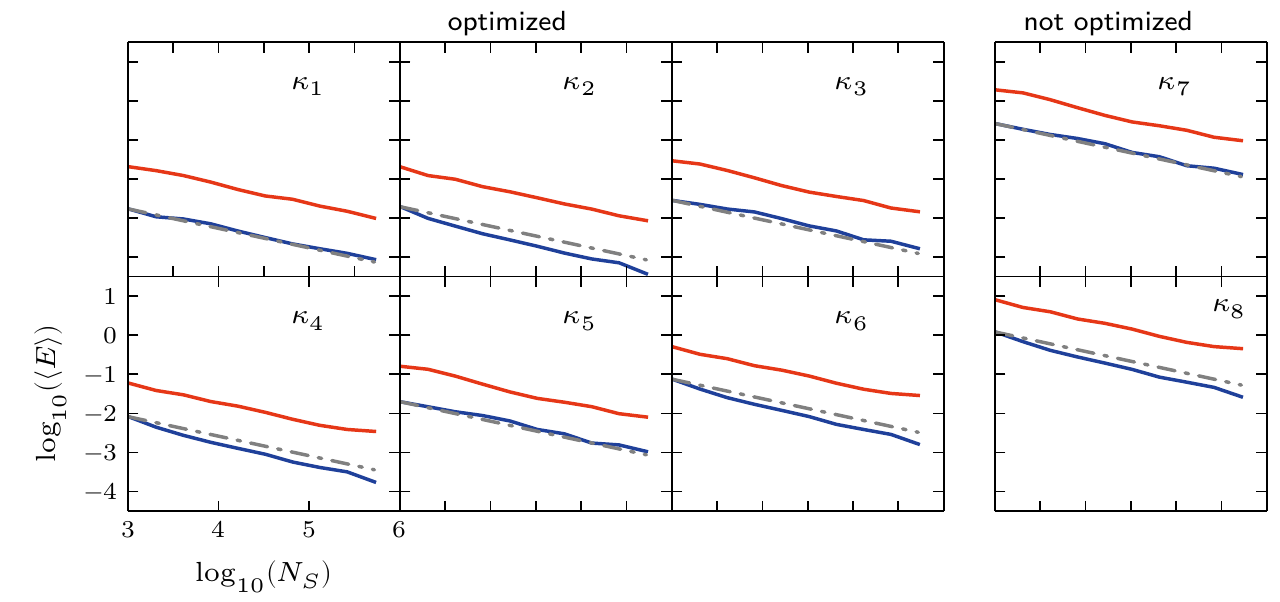}

\protect\caption{\label{fig:sde:nonlinear:combined:diffs-cumulants}Nonlinear SDE case,
showing the difference of the first $8$ cumulants from the steady-state
values at $t=25$ for the POS integrator using the combined method
(solid blue lines) and a reference explicit Euler integrator (dashed
red lines). Both integrators use $N_{T}=12.5\times10^{3}$ time steps.
The dash-dotted grey lines shows the slope for $\propto N_{S}^{-1/2}$. }
\end{figure}

In Fig.~\ref{fig:sde:nonlinear:combined:diffs-cumulants} we show
results for the cumulants in the nonlinear case, for completeness.
Unlike the linear case, where the improvement was confined to only
the explicitly optimized cumulants, we see that in the nonlinear case
even the non-optimized cumulants are optimized as well. This is because
the nonlinear equations couple all cumulants with each other. The
results show that POS results in dependency slightly better than $N_{S}^{-1/2}$
for the even cumulants. This is an interesting and noteworthy results
though we don't have an explanation for it yet.

In a separate analysis, we have integrated Eq.~(\ref{eq:nonlin_SDE})
until $t=4$ with a time-step of $\Delta t=4\times10^{-4}$ using
both the combined and the individual POS-method as well as an explicit
Ito-Euler reference method. We have compared the values of the first
$8$ consecutive moments and cumulants with that of a regular SDE
integrator with $N_{T}=8\times10^{4}$ time steps and $N_{S}=10^{9}$
trajectories, using a central difference integration method~\cite{Werner1997}.
The large number of trajectories for the regular SDE integrator ensured
that we obtained a sufficiently accurate estimate for the true value
of the moments and cumulants at $t=4$. The comparison showed that
individual POS method performed very similarly to the combined POS
method, resulting in a similarly good accuracy improvement. In the
interest of space, we are limiting our results to the combined POS
method from now on.

\subsection{Irregular drift, additive noise\label{sub:Irregular-drift,-additive}}

We now investigate the properties of an SDE whose functional behaviour
is not a regular polynomial. We chose the SDE

\begin{equation}
\mathrm{d}x=x\left(1-\left|x\right|\right)\mathrm{d}t+\mathrm{d}w,\label{eq:nonregular_SDE}
\end{equation}
where $\langle\mathrm{d}w\rangle=0$ and $\langle\mathrm{d}w^{2}\rangle=1$.

We can obtain the steady-state solutions for a given observable $\left\langle f\left(x\right)\right\rangle $
via

\begin{equation}
\left\langle f\left(x\right)\right\rangle _{ss}=\frac{\int_{-\infty}^{\infty}f\left(x\right)e^{x^{2}-\frac{2}{3}\left|x\right|x^{2}}\mathrm{d}x}{\int_{-\infty}^{\infty}e^{x^{2}-\frac{2}{3}\left|x\right|x^{2}}\mathrm{d}x}.
\end{equation}

\begin{figure}
\includegraphics{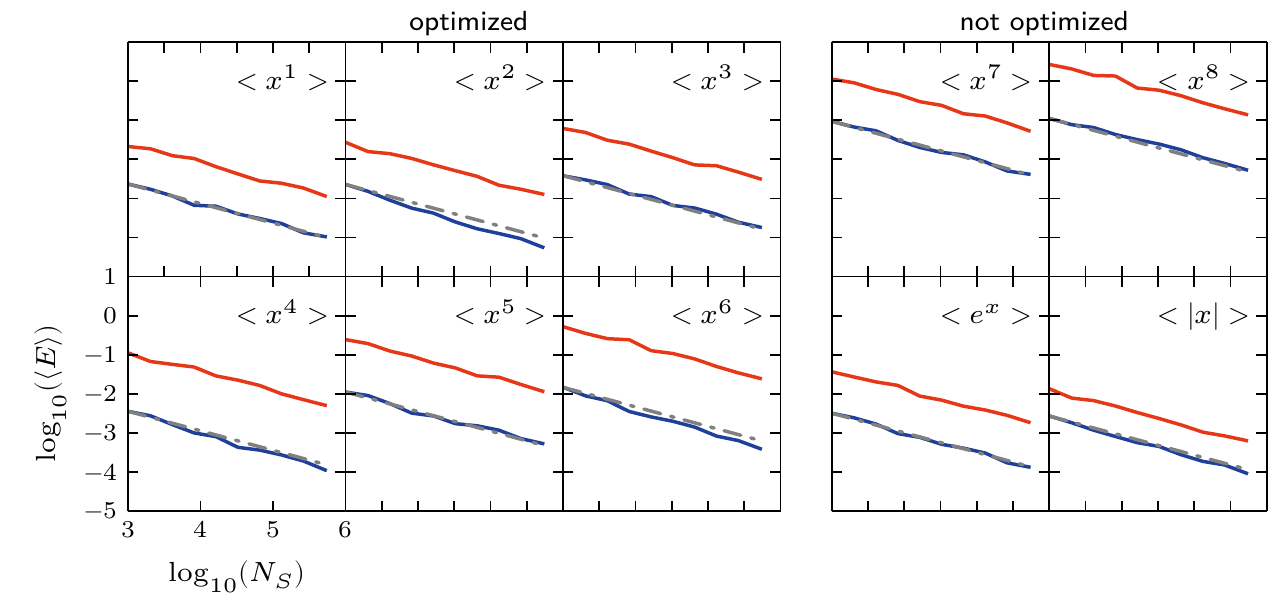}

\protect\caption{\label{fig:sde:nonregular:combined:diffs-moments}Nonlinear SDE case,
showing the difference of the first $8$ moments, the exponential
and absolute value function from the steady-state solution at $t=25$
for the POS integrator using the combined method (solid blue lines)
and a reference explicit Euler integrator (solid red lines). Both
integrators use $N_{T}=12.5\times10^{3}$ time steps. The dash-dotted
grey line shows the slope for $\propto N_{S}^{-1/2}$. }
\end{figure}

We run a simulation identical to the one in Section~\ref{sub:Nonlinear-drift,-additive-noise},
where we replace the SDE by Eq.~(\ref{eq:nonregular_SDE}). An analysis
analogous to that in Section~\ref{sub:Nonlinear-drift,-additive-noise}
reveals an equally good convergence to the steady-state values at
$t=25$. The results given in Fig.~\ref{fig:sde:nonregular:combined:diffs-moments}
show that the overall performance of POS for the non-regular SDE is
very similar to that for the nonlinear SDE studied in Section~\ref{sub:Nonlinear-drift,-additive-noise}.
Despite the fact that Eq.~(\ref{eq:nonregular_SDE}) is non-regular,
we observe virtually the same rate of improvement through POS as we
did for the regular SDE.

\begin{figure}
\includegraphics{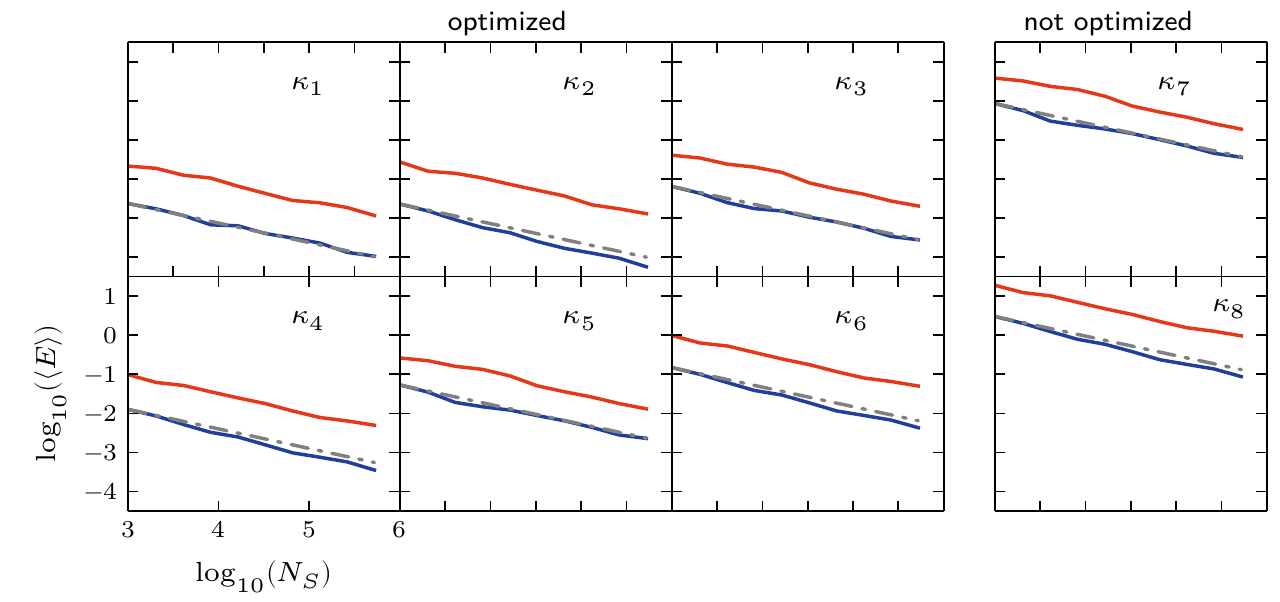}

\protect\caption{\label{fig:sde:nonregular:individual:diffs-cumulants}Nonlinear SDE
case, showing the difference of the first $8$ cumulants from the
accurate solution at $t=25$ for the POS integrator using the individual
method (solid blue lines) and a reference explicit Euler integrator
(solid red lines). Both integrators use $N_{T}=12.5\times10^{3}$
time steps. The dash-dotted grey line shows the slope for $\propto N_{S}^{-1/2}$
. }
\end{figure}

\section{Two-dimensional case\label{sec:Two-dimensional-case}}

Finally, to demonstrate that POS can be applied to SDE's with dimension
higher than $1$, we briefly consider the so-called laser equation.
This is a simplified model for the single-mode quantum statistics
in a laser system. Although it can be written as a single complex
equation, it is an example of an SDE with dimension higher than one
in real variables. We consider
\begin{equation}
da=\left(1-\left|a\right|^{2}\right)a\ dt+b\ dW_{c},\label{eq:laser_equation}
\end{equation}
where both $a$ and $dW_{c}$ are complex-valued. The noise obeys
$\left\langle dW_{c}dW_{c}^{*}\right\rangle =2dt$. Here $a$ is proportional
to the mode amplitude, and $b$ is a real-valued parameter that depends
on the steady-state photon number of the laser mode~\cite{Drummond2014-quantum}.
By treating the real and imaginary parts of $a$ as separate variables,
Eq.~(\ref{eq:laser_equation}) constitutes a two-dimensional real-valued
SDE.

The scaled photon number $n$ is given by $n=\left|a\right|^{2}$.
A simple calculation reveals the steady-state value for $n$ to be
\begin{equation}
n_{ss}=1+\sqrt{\frac{2}{\pi}}\frac{b}{\exp\left[1/\left(2b^{2}\right)\right]\left[1+\text{erf}\left(1/\left(\sqrt{2}b\right)\right)\right]}\,.
\end{equation}
Using $N_{S}=131072$ trajectories for the real and imaginary part
$\Re\left\{ a\right\} $ and $\Im\left\{ a\right\} $, resulting in
$N_{S}=262144$ trajectories in total, we initialize the ensemble
with normally-distributed numbers with mean $0$ and standard deviation
$\frac{1}{\sqrt{2}}$ for $\Re\left\{ a\right\} $ and $\Im\left\{ a\right\} $
and carry out a POS-optimized integration using the combined method
as well as an explicit Ito-Euler reference integration with an integration
time of $t=25$. An analysis similar to that in Sections~\ref{sub:Nonlinear-drift,-additive-noise}
and~\ref{sub:Irregular-drift,-additive} reveals an equally good
convergence to the steady-state at $t=25$. We then compare the value
of $n=\left|a\right|^{2}$ at $t=25$ with the steady-state value.

\subsection{Error in the equilibrium photon number}

We are optimizing (constraining) every moment and cross-moment for
the real and imaginary part of $a$ up to order 4. In other words,
we are optimizing the first 4 consecutive moments of $\Re\left\{ a\right\} $,
of $\Im\left\{ a\right\} $ and cross-moments of the form $\Re\left\{ a\right\} ^{n}\Im\left\{ a\right\} ^{m}$,
with a combined power $n+m$ not greater than 4. This results in a
total of 18 constrained moments. The POS optimization includes the
optimization of the initial trajectories.

We vary the parameter $b$ over a range from $1\cdot10^{-2}$ to $10.24$
in $11$ steps. As in Sections~\ref{sub:Nonlinear-drift,-additive-noise}
and~\ref{sub:Irregular-drift,-additive} we average over $120$ independent
test runs for the POS integration and $600$ runs for the reference
integration.

\begin{figure}
\centering{}\includegraphics[width=0.5\columnwidth]{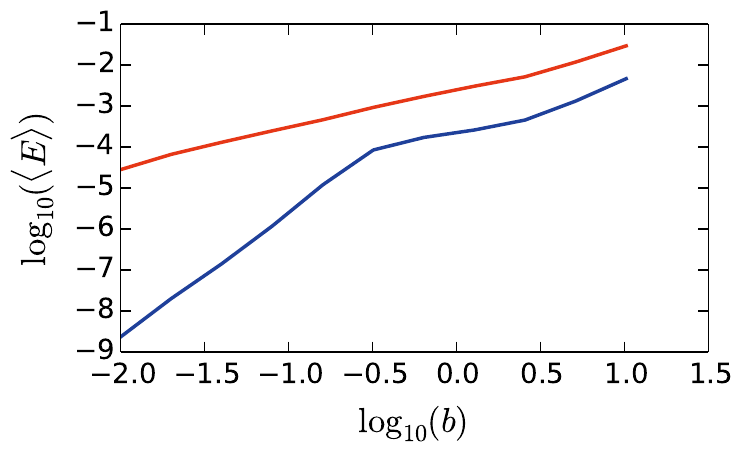}\protect\caption{\label{fig:sde:laser equation}Accuracy of the scaled photon number
$n_{s}$ for the laser equation using a reference explicit Euler integrator
(red line) and a POS-optimized integration (blue line) as a function
of the parameter $b$. Both integrators use $N_{T}=12.5\times10^{3}$
time steps.}
\end{figure}

As a typical physical example of a scientifically important quantity
in a stochastic calculation, Fig.~\ref{fig:sde:laser equation} shows
the error in calculating $n_{s}$, the equilibrium photon number in
the laser cavity. Here we use the combined Ito POS method (blue line),
and a non-optimized stochastic integration (red line) for the laser
equation for different values of $b$.

This example shows that POS can be successfully used to improve the
accuracy of complex-valued SDEs and, more generally speaking, real-valued
multivariate SDEs by several orders of magnitude.

\section{Conclusions\label{sec:conclusions}}

In summary, we have proposed and implemented a novel variance reduction
technique for solving stochastic differential equations, using parallel
optimized sampling. The essential feature is that it unifies moment
hierarchy and independent stochastic methods.

A finite moment hierarchy condition is imposed as a nonlinear constraint
on the random noises generated at each step in the integration. This
gives a dramatic reduction in sampling error for all moments calculated.
In the case of linear equations, we find that the low order moments
have sampling errors reduced to machine accuracy. While higher order
moments also have their errors reduced, the higher order cumulants
are not affected.

For nonlinear equations, the error reduction is not as large, but
it occurs over all moments and cumulants studied, even including the
non-optimized cumulants of higher order than the optimization limit.
We emphasize that the proposed algorithms have a general applicability
to all types of stochastic equations, and can be extended in principle
to higher order methods as well as the simple Euler integration treated
here.

\section*{Acknowledgements}

We wish to acknowledge financial support from the Australian Research
Council through the Discovery Grant program, as well as useful discussions
with Rodney Polkinghorne, Piotr Deuar and Brian Dalton.

\bibliographystyle{elsarticle-num}
\bibliography{POS}

\end{document}